\input ./style/arxiv-general.cfg
\documentclass[aos,MSNbibl,dvips]{arximspdf}
\makeatletter
   \@ifpackageloaded{graphicx}{}{\usepackage{graphicx}}
\makeatother
\usepackage{prettyref}

\doi{10.1214/15-AOS1332}
\volume{43}
\issue{5}
\pubyear{2015}
\firstpage{2168}
\lastpage{2197}
\docsubty{FLA}

\makeatletter

\newcommand{\rrvert}{\vert}

\newcommand{\rrVert}{\Vert}
\newcommand{\llvert}{\vert}
\newcommand{\llVert}{\Vert}
\newcommand{\eqref}[1]{(\ref{#1})}

\newcommand{\reals}{{\mathbb{R}}}
\newcommand{\supp}{\operatorname{supp}}
\newcommand{\eexp}{{e}}

\newcommand{\Tr}{\operatorname{\mathsf{Tr}}}
\newcommand{\diag}{\mathop{\operatorname{diag}}}

\newcommand{\Expect}{\mathbb{E}}
\newcommand{\Cov}{\operatorname{Cov}}

\newcommand{\calF}{{\mathcal{F}}}

\newcommand{\eps}{\varepsilon}
\newcommand{\wh}{\widehat}
\newcommand{\wt}{\widetilde}
\newcommand{\Var}{\operatorname{Var}}
\newtheorem{theorem}{Theorem}

\newtheorem{lemma}{Lemma}

\newproclaim{definition}{Definition}
\newproclaim{condition}{Condition}
\newproclaim{assumption}{Assumption}
\newproclaim{question}{Question}
\newproclaim{remark}{Remark}
\newproclaim{example}{Example}
\makeatother

\begin{document}
\begin{frontmatter}

\title{Minimax estimation in sparse canonical correlation~analysis}
\runtitle{Sparse CCA}

\begin{aug}
\author[A]{\fnms{Chao}~\snm{Gao}\thanksref{T1,m1}\ead[label=e1]{chao.gao@yale.edu}},
\author[B]{\fnms{Zongming}~\snm{Ma}\thanksref{T2,m2}\ead[label=e2]{zongming@wharton.upenn.edu}\ead[label=u2,url]{http://www-stat.wharton.upenn.edu/\textasciitilde zongming}},
\author[C]{\fnms{Zhao}~\snm{Ren}\thanksref{T1,m3}\ead[label=e3]{zren@pitt.edu}\ead[label=u3,url]{http://www.pitt.edu/\textasciitilde zren}}\\
\and
\author[A]{\fnms{Harrison H.}~\snm{Zhou}\corref{}\thanksref{T1,m1}\ead[label=e4]{huibin.zhou@yale.edu}\ead[label=u4,url]{http://www.stat.yale.edu/\textasciitilde hz68}}
\thankstext{T1}{Supported in part by NSF Career Award
DMS-0645676 and NSF FRG Grant DMS-08-54975.}
\thankstext{T2}{Supported in part by NSF Career
Award DMS-1352060 and the Claude Marion Endowed
Faculty Scholar Award of the Wharton School.}

\runauthor{Gao, Ma, Ren and Zhou}
\affiliation{Yale University\thanksmark{m1}, University of
Pennsylvania\thanksmark{m2} and University of Pittsburgh\thanksmark{m3}}
\address[A]{C. Gao\\
H. H. Zhou\\
Department of Statistics\\
Yale University\\
New Haven, Connecticut 06511\\
USA\\
\printead{e1}\\
\phantom{E-mail:\ }\printead*{e4}\\
\printead{u4}}
\address[B]{Z. Ma\\
Department of Statistics\\
The Wharton School\\
University of Pennsylvania\\
Philadelphia, Pennsylvania 19104\\
USA\\
\printead{e2}\\
\printead{u2}}
\address[C]{H. Ren\\
Department of Statistics\\
University of Pittsburgh\\
Pittsburgh, Pennsylvania 15260\\
USA\\
\printead{e3}\\
\printead{u3}}
\end{aug}

%
\received{\smonth{5} \syear{2014}}
%
\revised{\smonth{2} \syear{2015}}

%
\begin{abstract}
Canonical correlation analysis is a widely used multivariate
statistical technique
for exploring the relation between two sets of variables.
This paper considers the problem of estimating the leading
canonical correlation directions in high-dimensional settings.
Recently, under the
assumption that the leading canonical correlation directions are
sparse, various
procedures have been proposed for many high-dimensional applications
involving massive data sets.
However, there has been few theoretical justification available in
the literature. In this paper, we establish rate-optimal
nonasymptotic minimax estimation
with respect to an appropriate loss function for a wide range of model spaces.
Two interesting phenomena are observed. First,
the minimax rates are not affected by the presence of nuisance parameters,
namely the covariance matrices of the two sets of random variables,
though they need to be estimated in the canonical correlation analysis problem.
Second, we allow the presence of the residual canonical correlation directions.
However, they do not influence the minimax rates under a mild condition
on eigengap.
A generalized sin-theta theorem and an empirical process bound for
Gaussian quadratic forms under rank constraint are used to establish
the minimax upper bounds,
which may be of independent interest.
\end{abstract}

%
\begin{keyword}[class=AMS]
\kwd[Primary ]{62H12}
\kwd[; secondary ]{62C20}
\end{keyword}
\begin{keyword}
\kwd{Covariance matrix}
\kwd{minimax rates}
\kwd{model selection}
\kwd{nuisance parameter}
\kwd{sin-theta theorem}
\kwd{sparse CCA (SCCA)}
\end{keyword}
%
\end{frontmatter}

\section{Introduction}
\label{sec:intro}
Canonical correlation analysis (CCA) \cite{hotelling36} is one of the
most classical and important tools in multivariate statistics \cite
{Anderson03,mkb}.
It has been widely used in various fields to explore the relation
between two sets of variables measured on the same sample.

On the population level,
given two random vectors $X\in\reals^p$ and $Y\in\reals^m$, CCA first
seeks two vectors $u_1\in\reals^p$ and $v_1\in\reals^m$ such that the
correlation between the projected variables $u_1'X$ and $v_1'Y$ is maximized.
More specifically, $(u_1, v_1)$ is the solution to the following
optimization problem:
%
\begin{equation}
\label{eq:cca-def} \max_{u\in\reals^p, v\in\reals^m} \Cov\bigl(u'X,
v'Y\bigr),\qquad \mbox{subject to}\quad \Var\bigl(u'X\bigr) =
\Var\bigl(v'Y\bigr) = 1,
\end{equation}
which is uniquely determined up to a simultaneous sign change when
there is a positive eigengap.
Inductively, once $(u_i, v_i)$ is found, one can further obtain
$(u_{i+1}, v_{i+1})$ by solving the above optimization problem
repeatedly subject to the extra constraint that
\[
\Cov\bigl(u'X, u_j'X\bigr) = \Cov
\bigl(v'Y, v_j'Y\bigr) = 0\qquad \mbox{for }
j=1,\ldots, i.
\]
Throughout the paper, we call the $(u_i,v_i)$'s canonical correlation
directions.
It was shown by Hotelling \cite{hotelling36} that the $(\Sigma_x^{1/2}
u_i, \Sigma_y^{1/2}v_i)$'s are the successive singular vector pairs of
%
\begin{equation}
\label{eq:cca-svd} \Sigma_x^{-{1/2}} \Sigma_{xy}
\Sigma_y^{-{1/2}},
\end{equation}
where $\Sigma_x = \Cov(X), \Sigma_y = \Cov(Y)$ and $\Sigma_{xy} =
\Cov(X,Y)$.
When one is only given a random sample $\{(X_i, Y_i):i=1,\ldots, n\}$ of
size $n$, classical CCA estimates the canonical correlation directions
by performing singular value decomposition (SVD) on the sample
counterpart of \eqref{eq:cca-svd} first and then premultiply the
singular vectors by the inverse of square roots of the sample
covariance matrices.
For fixed dimensions $p$ and $m$, the estimators are well behaved when
the sample size is large~\cite{Anderson99}.

However, in contemporary datasets, we typically face the situation
where the ambient dimension in which we observe data is very high while
the sample size is small.
The dimensions $p$ and $m$ can be much larger than the sample size $n$.
For example, in cancer genomic studies, $X$ and $Y$ can be gene
expression and DNA methylation measurements, respectively, where the
dimensions $p$ and $m$ can be as large as tens of thousands while the
sample size $n$ is typically no larger than several hundreds \cite{cancer12}.
When applied to datasets of such nature, classical CCA faces at least
three key challenges.
First, the canonical correlation directions obtained through classical
CCA procedures involve all the variables measured on each subject, and
hence are difficult to interpret.
Second, due to the amount of noise that increases dramatically as the
ambient dimension grows, it is typically impossible to consistently
estimate even the leading canonical correlation directions without any
additional structural assumption \cite{Johnstone08,Bao14}.
Third, successive canonical correlation directions should be orthogonal
with respect to the population covariance matrices which are
notoriously hard to estimate in high-dimensional settings.
Indeed, it is not possible to obtain a substantially better estimator
than the sample covariance matrix \cite{MaWu13} which usually behaves
poorly \cite{Johnstone01}.
So, the estimation of such nuisance parameters further complicates the
problem of high-dimensional CCA.


Motivated by genomics, neuroimaging and other applications,
there have been growing interests in imposing sparsity assumptions on
the leading canonical correlation directions.
See, for example, \cite
{wiesel08,witten09,parkhomenko09,hardoon2011sparse,le2009sparse,waaijenborg2009sparse,avants2010dementia,Wang14}
for some
recent methodological developments and applications.
By seeking sparse canonical correlation directions, the estimated
$(u_i, v_i)$ vectors only involve a small number of variables, and
hence are easier to interpret.

Despite these recent methodological advances, theoretical understanding
about the sparse CCA problem is lacking.
It is unclear whether the sparse CCA algorithms proposed in the
literature have consistency or certain rates of convergence if the
population canonical correlation directions are indeed sparse.
To the best of our limited knowledge, the only theoretical work
available in the literature is \cite{chen13}.
In this paper, the authors gave a characterization for the sparse CCA
problem and considered an idealistic single canonical pair model where
$\Sigma_{xy}$, the covariance between $X$ and $Y$, was assumed to have
a rank one structure. They reparametrized $\Sigma_{xy}$ as follows:
%
\begin{equation}
\label{eq:scp-model} \Sigma_{xy} = \Sigma_x \lambda
uv' \Sigma_y,
\end{equation}
where $\lambda\in(0,1)$ and $u'\Sigma_x u = v'\Sigma_y v = 1$.
It can be shown that $(u,v)$ is the solution to \eqref{eq:cca-def}, so
that they are the leading canonical correlation directions.
It is worth noting that without knowledge of $\Sigma_x$ and $\Sigma_y$,
one is not able to obtain (resp., estimate) $(u,v)$ by simply applying
singular value decomposition to $\Sigma_{xy}$ (resp., sample covariance
$\wh\Sigma_{xy}$).
Under this model, Chen et al. \cite{chen13} studied the minimax lower
bound for estimating the individual vectors $u$ and $v$, and proposed
an iterative thresholding approach for estimating $u$ and $v$,
partially motivated by \cite{Ma11}.
However, their results depend on how well the nuisance parameters
$\Sigma_x$ and $\Sigma_y$ can be estimated, which to our surprise,
turns out to be unnecessary as shown in this paper.

\subsection{Main contributions}
The main objective of the current paper is to understand the
fundamental limits of the sparse CCA problem from a decision-theoretic
point of view. Such an investigation is not only interesting in its own
right, but will also inform the development and evaluation of practical
methodologies in the future. The model considered in this work is very general.
As shown in \cite{chen13}, $\Sigma_{xy}$ can be reparametrized as follows:
%
\begin{equation}
\label{eq:CCA} \Sigma_{xy} = \Sigma_x \bigl(U\Lambda
V'\bigr) \Sigma_y\qquad \mbox{with } U'
\Sigma_x U = V'\Sigma_y V =
I_{\bar{r}},
\end{equation}
where $\bar{r} = \min(p,m)$, $\Lambda= \diag(\lambda_1,\ldots,
\lambda
_{\bar{r}})$ and $1 > \lambda_1\geq\cdots\geq\lambda_{\bar{r}}
\geq0$.
Then the successive columns of $U$ and $V$ are the leading canonical
correlation directions.
Therefore, \eqref{eq:CCA} is the most general model for covariance
structure, and sparse CCA actually means the leading columns of $U$ and
$V$ are sparse.

We can split $U\Lambda V'$ as
%
\begin{equation}
U\Lambda V' = U_1 \Lambda_1
V_1' + U_2 \Lambda_2
V_2', \label{eq:splitCCA}
\end{equation}
where $\Lambda_1 = \diag(\lambda_1,\ldots,\lambda_r), \Lambda_2 =
\diag
(\lambda_{r+1},\ldots,\lambda_{\bar{r}})$, $U_1\in\mathbb
{R}^{p\times
r}$, $V_1\in\mathbb{R}^{m\times r}$, $U_2\in\reals^{p\times r_2}$ and
$V_2\in\reals^{m\times r_2}$ for $r_2 = \bar{r}-r$.
In what follows, we call $(U_1, V_1)$ the \emph{leading} and $(U_2,
V_2)$ the \emph{residual} canonical correlation directions.
Since our primary interest lies in $U_1$ and $V_1$, both the covariance
matrices $\Sigma_x$ and $\Sigma_y$ and the residual canonical
correlation directions $U_2$ and $V_2$ are nuisance parameters in our
problem. This model is more general than (\ref{eq:scp-model})
considered in \cite{chen13}. It captures the situation in real practice
where one is interested in recovering the first few sparse canonical
correlation directions while there might be additional directions in
the population structure.

To measure the performance of a procedure, we propose to estimate the
matrix $U_1 V_1'$ under the following loss function:
%
\begin{equation}
\label{eq:loss} L\bigl(U_1 V_1',
\wh{U_1 V_1'}\bigr) = \bigl\|U_1
V_1' - \wh{U_1 V_1'}
\bigr\|_{\mathrm{F}}^2.
\end{equation}
We choose this loss function for several reasons.
First, even when the $\lambda_i$'s are all distinct, $U_1$ and $V_1$
are only determined up to a simultaneous sign change of their columns.
In contrast, the matrix $U_1 V_1'$ is uniquely defined as long as
$\lambda_r > \lambda_{r+1}$.
Second, \eqref{eq:loss} is stronger than the squared projection error
loss. For any matrix $A$, let $P_A$ stand for the projection matrix
onto its column space. If the spectra of $\Sigma_x$ and $\Sigma_y$ are
both bounded away from zero and infinity, then, in view of Wedin's
sin-theta theorem \cite{Wedin72}, any upper bound on the loss function
\eqref{eq:loss} leads to an upper bound on the loss functions $\|
P_{U_1} - \wh{P}_{U_1}\|_{\mathrm{F}}^2$ and $\|P_{V_1} - \wh{P}_{V_1}\|
_{\mathrm{F}}^2$ for
estimating the column subspaces of $U_1$ and $V_1$, which have been
used in the related problem of sparse principal component analysis
\cite
{CMW13a,Vu13}.
Third, this loss function comes up naturally as the key component in
the Kullback--Leibler divergence calculation for a special class of
normal distributions where $\Sigma_x = I_p$, $\Sigma_y = I_m$ and
$\lambda_{r+1} = \cdots= \lambda_{\bar{r}} = 0$ in \eqref{eq:CCA}.

We use weak-$\ell_q$ balls to quantify sparsity. Let $\|{(U_1)_{j*}} \|
$ denote the $\ell_2$ norm of the $j$th row of $U_1$, and
let $\|{(U_1)_{(1)*}} \|\geq\cdots\geq\|{(U_1)_{(p)*}} \|$ be the
ordered row norms.
One way to characterize the sparsity in $U_1$ (and $V_1$) is to look at
its weak-$\ell_q$ radius for some $q\in[0,2)$,
%
\begin{equation}
\label{eq:weak-lq} \|{U_1} \|_{q,w} = \max_{j\in[p]}
j \bigl\|{(U_1)_{(j)*}} \bigr\|^q
\end{equation}
under the tradition that $0^q = 0$.
For instance, in the case of exact sparsity, that is, $q = 0$, $\|{U_1}
\|
_{0,w}$ counts the number of nonzero rows in $U_1$.
When $q \in(0,2)$, \eqref{eq:weak-lq}~quantifies the decay of the
ordered row norms of $U_1$, which is a form of approximate sparsity.
Then we define the parameter space $\mathcal{F}_q (s_u, s_v, p,m, r,
\lambda; \kappa, M)$, as the collection of all covariance matrices
\[
\Sigma= %
\left[\matrix{ \Sigma_x & \Sigma_{xy}
\cr
\Sigma_{yx} & \Sigma_y} \right] %
\]
with the CCA structure (\ref{eq:CCA}) and (\ref{eq:splitCCA}),
which satisfies:
\begin{enumerate}
\item$U_1\in\mathbb{R}^{p\times r}$ and $V_1\in\mathbb{R}^{m\times r}$
satisfying $\|{U_1} \|_{q,w}\leq s_u$ and $\|{V_1} \|_{q,w}\leq s_v$;
\item$\llVert \Sigma_x^l\rrVert _{\mathrm{op}}\vee\llVert \Sigma
_y^l\rrVert _{\mathrm{op}%
}\leq M$ for $l=\pm1$;
\item$1>\kappa\lambda\geq\lambda_1\geq\cdots\geq\lambda_r\geq
\lambda> 0$.
\end{enumerate}
Throughout the paper, we assume $\kappa\lambda\leq1-c_0$ for some
absolute constant $c_0\in(0,1)$.
The key parameters $s_u, s_v, p,m, r$ and $\lambda$ are allowed to
depend on the sample size $n$, while $\kappa, M> 1$ are treated as
absolute constants.
Compared with the single canonical pair model \eqref{eq:scp-model} in
\cite{chen13}, where $\operatorname{rank}(\Sigma_{xy})=1$, in this
paper, the
rank of $\Sigma_{xy}$ can be as high as $p$ or $m$ and $r$ is allowed
to grow. In addition, we do not need any structural assumption on
$\Sigma_x$ and $\Sigma_y$ except for condition 2 on the largest and
smallest eigenvalues, which implies that $\Sigma_x$ and $\Sigma_y$ are
invertible.

Suppose we observe i.i.d. pairs $(X_1,Y_1),\ldots,(X_n, Y_n)\sim
N_{p+m}(0,\Sigma)$.
For two sequences $\{a_n\}$ and $\{b_n\}$ of positive numbers, we write
$a_n\asymp b_n$ if for some absolute constant $C>1$, $1/C \leq a_n/b_n
\leq C$ for all $n$.
By the minimax lower and upper bound results in Section~\ref{sec:result},
under mild conditions, we obtain the following tight nonasymptotic
minimax rates for estimating the leading canonical directions when $q = 0$:
%
\begin{eqnarray}
\label{eq:rate-0} %
&& \inf_{\wh{U_1 V_1'}} \sup
_{\Sigma\in\mathcal
{F}_0(s_u,s_v,p,m,r,\lambda)} \Expect_\Sigma\bigl\|U_1
V_1' - \wh{U_1 V_1'}
\bigr\|_{\mathrm{F}}^2
\nonumber
\\[-8pt]
\\[-8pt]
\nonumber
&&\qquad\asymp \frac{1}{n\lambda^2}\biggl( r(s_u+s_v) +
s_u \log\frac{\eexp p}{s_u} + s_v \log\frac{\eexp m}{s_v}
\biggr). %
\end{eqnarray}
In Section~\ref{sec:result}, we give a precise statement of this result
and tight minimax rates for the case of approximate sparsity, that is,
$q\in(0,2)$.

The result (\ref{eq:rate-0}) provides a precise characterization of the
statistical fundamental limit of the sparse CCA problem.
It is worth noting that the conditions required for \eqref{eq:rate-0}
do not involve any additional assumptions on the nuisance parameters
$\Sigma_x, \Sigma_y, U_2$ and $V_2$.
Therefore, we are able to establish the remarkable fact that the
fundamental limit of the sparse CCA problem is \emph{not} affected by
those nuisance parameters. This optimality result can serve as an
important guideline to evaluate procedures proposed in the literature.

To obtain minimax upper bounds, we propose an estimator by optimizing
canonical correlation under sparsity constraints.
A key element in analyzing the risk behavior of the estimator is a
generalized sin-theta theorem. See Theorem~\ref{thm:sintheta} in
Section~\ref{sec:key}. The theorem is of interest in its own right and
can be useful in other problems where matrix perturbation analysis is needed.
It is worth noting that the proposed procedure does \emph{not} require
sample splitting, which was needed in \cite{CMW13a}. We bypass sample
splitting by establishing a new empirical process bound for the supreme
of Gaussian quadratic forms with rank constraint. See Lemma~\ref
{lem:EP} in Section~\ref{sec:key}.
The estimator is shown to be minimax rate optimal by establishing
matching minimax lower bounds based on a local metric entropy approach
\cite{Lecam73,Birge83,Yang99,CMW13a}.

\subsection{Connection to and difference from sparse PCA}
The current paper is related to the problem of sparse principal
component analysis (PCA), which has received a lot of recent attention
in the literature.
Most literature on sparse PCA considers the spiked covariance model
\cite{tb99jrssb,Johnstone01} where one observes an $n\times p$ data
matrix, each row of which is independently sampled from a normal
distribution $N_p(0,\Sigma_0)$ with
%
\begin{equation}
\label{eq:spiked-cov} \Sigma_0 = V\Lambda V' +
\sigma^2 I_p.
\end{equation}
Here, $V\in\reals^{p\times r}$ has orthonormal column vectors which
are assumed to be sparse and $\Lambda= \diag(\lambda_1,\ldots
,\lambda
_r)$ with $\lambda_1\geq\cdots\geq\lambda_r >0$.
Since the first $r$ eigenvalues of $\Sigma_0$ are $\{\lambda_i+\sigma
^2\}_{i=1}^r$ and the rest are all $\sigma^2$, the $\lambda_i$'s are
referred as ``spikes,'' and hence the name of the model.
Johnstone and Lu \cite{JohnstoneLu09} proposed a diagonal thresholding
estimator of the sparse principal eigenvector which is provably consistent
for a range of sparsity regimes.
For fixed $r$, Birnbaum et al. \cite{Birnbaum12} derived minimax rate
optimal estimators for individual sparse principal eigenvectors, and
Ma \cite{Ma11} proposed to directly estimate sparse principal
subspaces, that is, the span of $V$, and constructed an iterative
thresholding algorithm for this purpose which is shown to achieve near
optimal rate of convergence adaptively.
Cai et al. \cite{CMW13a} studied minimax rates and adaptive estimation
for sparse principal subspaces with little constraint on $r$.
See also \cite{Vu13} for the case of a more general model.
In addition, variable selection, rank detection, computational
complexity and posterior contraction rates of sparse PCA have been
studied. See, for instance, \cite{Amini09,CMW13,Berthet13,Gao13} and
the references therein.

Compared with sparse PCA, the sparse CCA problem studied in the current
paper is different and arguably more challenging in three important ways.
\begin{itemize}
\item In sparse PCA, the sparse vectors of interest, that is, the
columns of $V$ in \eqref{eq:spiked-cov} are normalized with respect to
the identity matrix.
In contrast, in sparse CCA, the sparse vectors of interest, that is,
the columns of $U$ and $V$ are normalized with respect to $\Sigma_x$
and $\Sigma_y$, respectively, which are not only unknown but also hard
to estimate in high-dimensional settings. The necessity of
normalization with respect to nuisance parameters adds on to the
difficulty of the sparse CCA problem.

\item In sparse PCA, especially in the spiked covariance model, there
is a clean separation between ``signal'' and ``noise'': the signal is
in the spiked part and the rest are noise.
However, in the parameter space considered in this paper, we allow the
presence of residual canonical correlations $U_2 \Lambda_2 V_2'$, which
is motivated by the situation statisticians face in practice.
It is highly nontrivial to show that the presence of the residual
canonical correlations does not influence the minimax estimation rates.

\item The covariance structures in sparse PCA and sparse CCA have both
sparsity and low-rank structures. However, there is a subtle difference
between the two. In sparse PCA, the sparsity and orthogonality of $V$
in (\ref{eq:spiked-cov}) are coherent. This means that the columns of
$V$ are sparse and orthogonal to each other simultaneously. Such
convenience is absent in the sparse CCA problem. It is implied from
(\ref{eq:CCA}) that $\Sigma_x^{1/2}U_1$ and $\Sigma_y^{1/2}V_1$ have
orthogonal columns, while it is the columns of $U_1$ and $V_1$ that are
sparse. The orthogonal columns and the sparse columns are different.
The consequence is that in order to estimate the sparse matrices $U_1$
and $V_1$, we must appeal to the orthogonality in the nonsparse
matrices $\Sigma_x^{1/2}U_1$ and $\Sigma_y^{1/2}V_1$, even when the
matrices $\Sigma_x$ and $\Sigma_y$ are unknown. If we naively treat
sparse CCA as sparse PCA, the procedure can be inconsistent (see the
simulation results in \cite{chen13}).
\end{itemize}

\subsection{Organization of the paper}
The rest of the paper is organized as follows. Section~\ref{sec:result}
presents the main results of the paper, including upper bounds in
Section~\ref{sec:upper} and lower bounds in Section~\ref{sec:lower}.
Section~\ref{sec:discussion} discusses some related issues.
The proofs of the minimax upper bounds are gathered in Section~\ref
{sec:proof}, with some auxiliary results and technical lemmas proved in
Section~\ref{sec:aux-proof}. The proof of the lower bounds and some
further technical lemmas are given in the supplementary material~\cite{supp2}.


\subsection{Notation}

For any matrix $A = (a_{ij})$,
the $i$th row of $A$ is denoted by ${A}_{{i} *}$ and the $j$th column by
${A}_{* {j}}$.
For a positive integer $p$, $[p]$ denotes the index set $\{1, 2, \ldots,
p\}$. For any set $I$, $|I|$ denotes its cardinality and ${I^{\mathrm{
c}}}$ its
complement.
For two subsets $I$ and $J$ of indices, we write $A_{IJ}$ for the
$|I|\times|J|$ submatrices formed by $a_{ij}$ with $(i,j) \in I \times J$.
When $I$ or $J$ is the whole set, we abbreviate it with an $*$, and so
if $A\in\reals^{p\times k}$, then ${A}_{{I} *} = A_{I [k]}$ and
${A}_{* {J}} = A_{[p] J}$.
For any square matrix $A = (a_{ij})$, denote its trace by $\Tr(A) =
\sum_{i}a_{ii}$.
Moreover, let $O(p,k)$ denote the set of all $p\times k$ orthonormal
matrices and $O(k)=O(k,k)$.
For any matrix $A \in\reals^{p \times k}$, $\sigma_i(A)$ stands for
its $i$th largest singular value. The Frobenius norm and the operator
norm of $A$ are defined as $\|{A} \|_{{\mathrm{ F}}}=\sqrt{\Tr
(A'A)}$ and $\|{A} \|_{{\mathrm{ op}}}=\sigma_1(A)$, respectively.
The support of $A$ is defined as
$\supp
(A)=\{i\in[n]: \|A_{i*}\|>0\}$. The trace inner product of two matrices
$A,B\in\reals^{p \times k}$ is defined as $ \langle A, B
\rangle=\Tr(A'B)$.
For any number $a$, we use ${\lceil{a} \rceil}$ to denote the
smallest integer
that is no smaller than $a$.
For any two numbers $a$ and $b$, let $a\vee b = \max(a,b)$ and
$a\wedge
b = \min(a,b)$.
For any event $E$, we use ${\mathbf{1}_{\{{E}\}}}$ to denote its
indicator function.
We use $\mathbb{P}_{\Sigma}$ to denote the probability distribution of
$N_{p+m}(0,\Sigma)$ and $\mathbb{E}_{\Sigma}$ for the associated expectation.

\section{Main results}
\label{sec:result}

In this section, we state the main results of the paper. In Section~\ref{sec:upper}, we introduce a method to estimate the leading canonical
correlation directions. Minimax upper bounds are obtained. Section~\ref
{sec:lower} gives minimax lower bounds which match the upper bounds up
to a constant factor. We abbreviate the parameter space $\mathcal
{F}_q(s_u,s_v,p,m,r,\lambda;\kappa,M)$ as $\mathcal{F}_q$.

\subsection{Upper bounds} \label{sec:upper}

The main idea of the estimator proposed in this paper is to maximize
the canonical correlations under sparsity constraints. Note that the
SVD approach of the classical CCA \cite{hotelling36} can be written in
the following optimization form:
%
\begin{equation}
\max_{(A,B)} \Tr\bigl(A'\wh{\Sigma}_{xy}B
\bigr) \qquad\mbox{s.t.}\quad A'\wh{\Sigma }_xA=B'
\wh{\Sigma}_yB=I_r. \label{eq:optim}
\end{equation}
We generalize (\ref{eq:optim}) to the high-dimensional setting by
adding sparsity constraints.

Since the leading canonical correlation directions $(U_1,V_1)$ are weak
$\ell_q$ sparse, we introduce effective sparsity for $q\in[0,2)$,
which plays a key role in defining the procedure.
Define
%
\begin{eqnarray}
\label{eq:x-q-u} x_q^u & =&\max \biggl\{0\leq x\leq p: x
\leq s_u 
\biggl( \frac{n\lambda^2}{r+\log(ep/x)} \biggr)
^{q/2} \biggr\},
\\
\label{eq:x-q-v} x_q^v & =&\max
\biggl\{0\leq x\leq m: x\leq s_v 
\biggl( \frac{n\lambda^2}{r+\log(em/x)}
\biggr) 
^{q/2} \biggr\}.
\end{eqnarray}
The effective sparsity of $U_1$ and $V_1$ are defined as
%
\begin{equation}
k_q^{u}={\bigl\lceil{x_q^u} \bigr
\rceil}, \qquad k_q^v={\bigl\lceil{x_q^v}
\bigr\rceil}. \label{eq:effectivesparsity}
\end{equation}
For $j\geq k_q^u$, it can be shown that
\[
\bigl\Vert(U_1)_{(j)*}\bigr\Vert\leq \biggl(\frac{r+\log(ep/k_q^u)}{n\lambda
^2}
\biggr)^{1/2},
\]
for which the signal is not strong enough to be recovered from the
data. It holds similarly for $V_1$.

For $n$ i.i.d. observations $(X_i,Y_i)$, $i\in[n]$, we compute
the sample covariance matrix
\[
\widehat{\Sigma}=%
\left[\matrix{ \widehat{\Sigma}_x
& \widehat{\Sigma}_{xy}
\cr
\widehat{\Sigma}_{yx} & \widehat{
\Sigma}_{y}%
}\right] %
.
\]
The estimator $(\wh{U}_1,\wh{V}_1)$ for $(U_1,V_1)$, the leading $r$
canonical correlation directions, is defined as a solution to the
following optimization problem:
%
\begin{eqnarray}
\label{eq:pickset} %
 &&\max_{(A,B)}  \Tr
\bigl(A'\wh{\Sigma}_{xy}B\bigr)
\nonumber
\\[-8pt]
\\[-8pt]
\nonumber
&&\qquad\mbox{s.t.}\quad  A'\wh{\Sigma}_xA=B'\wh{
\Sigma}_yB=I_r\mbox{ and }\|{A} \|_{0,w}=
k_q^u, \|{B} \|_{0,w}= k_q^v.
\end{eqnarray}
When $q=0$, we have $k_q^u=s_u$ and $k_q^v=s_v$. Then the program (\ref
{eq:pickset}) is just a slight generalization of the classical approach
of \cite{hotelling36} with additional $\ell_0$ constraints $\|{A} \|
_{0,w}=s_u$ and $\|{B} \|_{0,w}=s_v$. By the definition of the
parameter space, it is also natural to impose the $\ell_q$ constraints
$\|{A} \|_{q,w}\leq s_u$ and $\|{B} \|_{q,w}\leq s_v$. Such constraints
were used by \cite{Vu13} to solve the sparse PCA problem. However,
their upper bounds require more assumptions due to the difficulty in
analyzing $\ell_q$ constraints. We use $\ell_0$ constraints on the
effective sparsity and obtain the optimal upper bound under minimal assumptions.

Set
%
\begin{equation}
\eps_n^{2}=\frac{1}{n\lambda^{2}} \biggl(r
\bigl(k_{q}^{u}+k_{q}^{v}
\bigr)+k_{q}^{u} 
\log\frac{ep}{k_q^u}+k_{q}^{v}
\log\frac{em}{k_q^v} \biggr), \label{eq:DefEpsilon}
\end{equation}
which is the minimax rate to be shown later.

\begin{theorem}
\label{thm:upperbound1} We assume that
%
\begin{eqnarray}
\varepsilon_n^2 &\leq&c, \label{eq:ass1}
\\
\lambda_{r+1} &\leq&c\lambda, \label{eq:ass2}
\end{eqnarray}
for some sufficiently small constant $c\in(0,1)$. For any constant
$C'>0$, there exists a constant $C>0$ only depending on $M,q,\kappa$
and $C'$,
such that for any $\Sigma\in\mathcal{F}_q$,
\[
\bigl\| \widehat{U}_{1}\widehat{V}_{1}^{\prime}-U_{1}V_{1}^{\prime
}
\bigr\| _{\mathrm{F}}^{2}\leq C\varepsilon_n^2,
\]
with $\mathbb{P}_{\Sigma}$-probability at least $1-\exp(-C^{\prime
}(k_{q}^{u}+\log
(ep/k_{q}^{u})))-\exp(-C^{\prime}(k_{q}^{v}+\log
(em/k_{q}^{v})))$.
\end{theorem}

\begin{remark}
It will be shown in Section~\ref{sec:lower} that  assumption (\ref
{eq:ass1}) is necessary for consistent estimation.  Assumption (\ref
{eq:ass2}) implies $\lambda_{r+1}\leq c\lambda_r$ for $c\in(0,1)$,
such that the eigengap is lower bounded as $\lambda_r-\lambda
_{r+1}\geq
(1-c)\lambda_r>0$.
\end{remark}

\begin{remark}
The upper bound $\varepsilon_n^2$ has two parts. The first part $\frac
{1}{n\lambda^2} (r(k_q^u+k_q^v) )$ is caused by low rank
structure, and the second part $\frac{1}{n\lambda^2} (k_q^u\log
(ep/k_q^u)+k_q^v\log(em/k_q^v) )$ is caused by sparsity. If $r\leq
\log(ep/k_q^u)\wedge\log(em/k_q^v)$, the second part dominates, while
the first part dominates if $r\geq\log(ep/k_q^u)\vee\log(em/k_q^v)$.
\end{remark}

\begin{remark}
The upper bound does not require any structural assumption on the
marginal covariance matrices $\Sigma_x$ and $\Sigma_y$ other than
bounds on the largest and the smallest eigenvalues.
Although in the high-dimensional setting, the sample covariance $\wh
{\Sigma}_x$ and $\wh{\Sigma}_y$ are not good estimators of the matrices
$\Sigma_x,\Sigma_y$, the normalization constraints $A'\wh{\Sigma
}_xA=B'\wh{\Sigma}_yB=I_r$,\vspace*{1pt} together with the sparsity of $A,B$, only
involve submatrices of $\wh{\Sigma}_x$ and $\wh{\Sigma}_y$. Under the
assumption (\ref{eq:ass1}), it can be shown that a $k_q^u\times k_q^u$
submatrix of $\wh{\Sigma}_x$ converges to\vspace*{-1pt} the corresponding submatrix
of $\Sigma_x$ with the rate $\sqrt{\frac{k_q^y\log(ep/k_q^u)}{n}}$
under operator norm uniformly over all $k_q^u\times k_q^u$ submatrices.
Similar results hold for $\wh{\Sigma}_y$ and $\Sigma_y$. See Lemma~\ref
{lem:covdeviation45} in Section~\ref{sec:bias-pf}. These rates are
dominated by the minimax rate $\varepsilon_n$ in (\ref{eq:DefEpsilon}).
\end{remark}

\begin{remark}
One of the major difficulties of sparse CCA is the presence of the
unknown $\Sigma_x$ and $\Sigma_y$. Suppose $\Sigma_x$ and $\Sigma_y$
are known, one may work with the transformed data $\{(\Sigma_x^{-1}X_i,
\Sigma_y^{-1}Y_i):i=1,\ldots,n\}$. The cross-covariance of the transformed
data is $\Sigma_x^{-1}\Sigma_{xy}\Sigma_y^{-1}=U\Lambda V'$, which
is a
sparse matrix. When $\operatorname{rank}(\Sigma_{xy})=1$, algorithms
such as
\cite{Yang11,chen13} can obtain the sparse singular vectors from
$\Sigma
_x^{-1}\wh{\Sigma}_{xy}\Sigma_y^{-1}$, which estimate $U_1$ and $V_1$
with optimal rate. When $\Sigma_x$ and $\Sigma_y$ are unknown,
structural assumptions are required on the covariance matrices in order
that $\Sigma_x^{-1}$ and $\Sigma_y^{-1}$ can be well estimated. Then
one can use the estimated $\Sigma_x^{-1}$ and $\Sigma_y^{-1}$ to
transform the data and apply the same sparse singular vector estimator
(see \cite{chen13}). However, unless $\Sigma_x=I_p$ and $\Sigma_y=I_m$,
this method cannot be extended to the case where $\operatorname
{rank}(\Sigma
_{xy})\geq2$, since the orthogonality of $U$ and $V$ is with respect
to general covariance matrices $\Sigma_x$ and $\Sigma_y$, respectively.
In the case where $\Sigma_x=I_p$ and $\Sigma_y=I_m$, the problem is
similar to sparse PCA, and the proof of Theorem~\ref{thm:upperbound1}
can be greatly simplified.
\end{remark}

To obtain the convergence rate in expectation, we propose a modified
estimator. The modification is inspired by the fact that $%
U_{1}V_{1}^{\prime}$ are bounded in Frobenius norm, because
%
\begin{equation}
\bigl\| U_{1}V_{1}^{\prime}\bigr\| _{\mathrm{F}}\leq\bigl\|
\Sigma_{x}^{-1/2}\bigr\| _{\mathrm{op}}\bigl\| \Sigma
_{x}^{1/2}U_{1}\bigr\| _{\mathrm{F}}\bigl\| \Sigma
_{y}^{1/2}V_{1}\bigr\| _{\mathrm{op}}\bigl\| \Sigma
_{y}^{-1/2}\bigr\| _{\mathrm{op}}\leq M\sqrt{r}. \label{eq:BoundTruth}
\end{equation}
Define $\widehat{U_{1}V_{1}^{\prime}}$ to be the
truncated version of $\widehat{U}_{1}\widehat{V}_{1}^{\prime}$ as
\[
\widehat{U_{1}V_{1}^{\prime}}=
\widehat{U}_{1}\widehat {V}_{1}^{\prime
}
{\mathbf{1}_{\{{\| \widehat{U}_{1}\widehat{V}_{1}^{\prime}\|
_{\mathrm{F}}\leq2M\sqrt{r}}\}}}. 
\]
The modification can be viewed as an improvement, because whenever\break $%
\| \widehat{U}_{1}\widehat{V}_{1}^{\prime}\| _{\mathrm{F}%
}>2M\sqrt{r}$, we have
\[
\bigl\| \widehat{U}_{1}\widehat{V}_{1}^{\prime}-U_{1}V_{1}^{\prime
}
\bigr\| _{\mathrm{F}}\geq\bigl\| \widehat{U}_{1}\widehat{V}%
_{1}^{\prime}
\bigr\| _{\mathrm{F}}-\bigl\| U_{1}V_{1}^{\prime
}\bigr\|
_{\mathrm{F}}\geq M\sqrt{r}\geq\bigl\|{0-U_1V_1'}
\bigr\|_{{\mathrm{ F}}}.
\]
Then it is better to estimate $U_{1}V_{1}^{\prime}$ by $0$.

\begin{theorem}
\label{thm:upperbound2}
Suppose (\ref{eq:ass1}) and (\ref%
{eq:ass2}) hold. In addition, assume that
%
\begin{eqnarray}
\exp\bigl(C_{1}\bigl(k_{q}^{u}+\log
\bigl(ep/k_{q}^{u}\bigr)\bigr)\bigr) &>&n
\lambda^{2}, \label{eq:TheoremExpAssup1}
\\
\exp\bigl(C_{1}\bigl(k_{q}^{v}+\log
\bigl(em/k_{q}^{v}\bigr)\bigr)\bigr) %
&>&n
\lambda^{2}, \label{eq:TheoremExpAssup2}
\end{eqnarray}
for some $C_{1}>0$, then there exists $C_{2}>0$ only depending on
$M,q,\kappa$ and $C_1$, such that
\[
\sup_{\Sigma\in\mathcal{F}_q}\mathbb{E}_{\Sigma}\bigl\| \widehat{U_{1}V_{1}^{\prime}}-U_{1}V_{1}^{\prime}
\bigr\| _{\mathrm{F}%
}^{2}\leq C_2\varepsilon_n^2.
\]
\end{theorem}

\begin{remark}
The assumptions (\ref{eq:TheoremExpAssup1}) and (\ref
{eq:TheoremExpAssup2}) imply the tail probability in Theorem~\ref
{thm:upperbound1} is sufficiently small. Once there exists a small
constant $\delta>0$, such that
\[
p\vee e^{k_q^u}\geq n^{\delta}\quad \mbox{and}\quad m\vee e^{k_q^v}
\geq n^{\delta}
\]
hold, then (\ref{eq:TheoremExpAssup1}) and (\ref{eq:TheoremExpAssup2})
also hold with some $C_1>0$. Notice that $p>n^{\delta}$ is commonly
assumed in high-dimensional statistics to have convergence results in
expectation. The assumption here is weaker than that.
\end{remark}

\subsection{Lower bounds} \label{sec:lower}

Theorems \ref{thm:upperbound1} and \ref{thm:upperbound2} show
that the procedure proposed in~(\ref{eq:pickset}) attains the rate
$\varepsilon_n^2$.
In this section, we show that this rate is optimal among all
estimators. More specifically, we show that the following minimax lower
bounds hold for $q\in[0,2)$.


\begin{theorem}
\label{thm:lower-bd-q}
Assume that
$1\leq r \leq\frac{k_q^u\wedge k_q^v}{2}$,
and that
%
\begin{equation}
n\lambda^2 \geq C_0 \biggl( r+\log{\eexp p\over k_q^u}
\vee\log{\eexp
m\over k_q^v} \biggr)\label{eq:assumptionlower}
\end{equation}
for some sufficiently large constant $C_0$.
Then there exists a constant $c>0$ depending only on $q$ and an
absolute constant $c_0$ such that the minimax risk for estimating
$U_1V_1'$ satisfies
\[
\inf_{(\wh{U}_1, \wh{V}_1)}\sup_{\Sigma\in\mathcal{F}_q} \Expect
_{\Sigma
}\bigl\|\wh{U}_1\wh{V}'_1 -
U_1 V_1'\bigr\|_{\mathrm{F}}^2
\geq c\varepsilon_n^2\wedge c_0.
\]
\end{theorem}

The proof of Theorem~\ref{thm:lower-bd-q} is given in the supplementary
material \cite{supp2}.

\begin{remark}
Assumption (\ref{eq:assumptionlower}) is necessary for consistent
estimation.
\end{remark}

\section{Discussion}
\label{sec:discussion}
We include below discussions on two related issues.

\subsection{Minimax rates for individual sparsity}

In this paper, we have derived tight minimax estimation rates for the
leading sparse canonical correlation directions where the sparsity is
depicted by the rapid decay of the ordered row norms in $U_1$ and $V_1$
(as characterized by the weak-$\ell_q$ notion).

Another interesting case of sparsity is when the individual column
vectors of $U_1$ and $V_1$ are sparse.
For instance, when
%
\begin{equation}
\label{eq:col-sparse} \|u_i\|_{q,w}\leq t_u\quad
\mbox{and}\quad \|v_i\|_{q,w}\leq t_v \qquad\forall i
\in[r],
\end{equation}
where the $\|\cdot\|_{q,w}$ is defined as in \eqref{eq:weak-lq} by
treating any $p$-vector as a $p\times1$ matrix.
Let $\calF_q^c = \calF_q^c(t_u, t_v, p, m, r, \lambda; \kappa, M)$ be
defined as in Section~\ref{sec:intro} following \eqref{eq:weak-lq} but
with the sparsity notion changed to that in \eqref{eq:col-sparse}.
Similar to \eqref{eq:x-q-u}--\eqref{eq:effectivesparsity}, let
\[
y_q^u = \max\biggl\{ 0\leq y\leq p: y\leq
t_u \biggl( \frac{n\lambda ^2}{\log
 (ep/(ry) )} \biggr)^{q/2} \biggr\},\qquad
j_q^u = {\bigl\lceil{y_q^u} \bigr
\rceil},
\]
and $y_q^v$ and $j_q^v$ be analogously defined.
Then we have:

\begin{theorem}
\label{thm:col-sparse}
Assume that $1\leq r\leq\frac{j_q^u\wedge j_q^v}{2}$, $2r j_q^u\leq
p$, $2r j_q^v\leq m$ and $n\lambda^2 \geq C_0(r+ \log\frac{ep}{rj_q^u}
\vee\log\frac{em}{rj_q^v})$ for some sufficiently large constant $C_0$.
Then there is a constant $c>0$ depending only on $q$ and an absolute
constant $c_0>0$ such that
%
\begin{equation}\qquad
\label{eq:col-lowbd} \inf_{\wh{U_1 V_1'}} \sup_{\Sigma\in\calF_q^c}
\Expect_\Sigma\bigl\| U_1 V_1' -
\wh{U_1 V_1'}\bigr\|_{\mathrm{F}}^2
\geq c_0 \wedge\frac{c}{n\lambda^2}r\biggl( j_q^u
\log\frac{\eexp p}{r
j_q^u } + j_q^v \log\frac{\eexp m}{r j_q^v }
\biggr).
\end{equation}
If in addition $rj_q^u\leq p^{1-\alpha}$, $r j_q^v \leq m^{1-\alpha}$
for some small $\alpha\in(0,1)$, $r \leq C \log(p\wedge m)$ for some
$C > 0$ and the conditions of Theorem~\ref{thm:upperbound2} are satisfied
with $k_q^u = r j_q^u$ and $k_q^v = r j_q^v$,
then a matching upper bound is achieved by the estimator in
Theorem~\ref{thm:upperbound2}
with $k_q^u = r j_q^u$ and $k_q^v = r j_q^v$.
\end{theorem}

The proof of Theorem~\ref{thm:col-sparse} is given in the supplementary
material \cite{supp2}.
The lower bound (\ref{eq:col-lowbd}) for individual sparsity is larger
than the minimax rate (\ref{eq:DefEpsilon}) for joint sparsity when
$t_u=s_u$ and $t_v=s_v$.

\subsection{Adaptation, computation and some recent work}

The main purpose of proposing the estimator in \eqref{eq:pickset} is to
determine the minimax estimation rates in sparse CCA problem
under weak assumptions.
Admittedly, it requires the knowledge of parameter space and is
computationally intensive.

Designing adaptive and computationally efficient procedures to achieve
statistically optimal performance is an interesting and important
research direction.
Built upon the insights developed in the current paper, Gao et al.
\cite
{gao14b} have proposed an adaptive and efficient procedure for sparse CCA.
The procedure first obtains a crude estimator via a convex relaxation
of the problem \eqref{eq:pickset} here which is then refined by a group
sparse linear regression.
The resulting estimator achieves optimal rates of convergence in
estimating the leading sparse canonical directions under a prediction
loss without imposing any structural assumption on $\Sigma_x$ and
$\Sigma_y$, when the residual directions are absent.
Notably, the procedure in \cite{gao14b} requires a larger sample size
than in the present paper, which has been shown to be essentially
necessary for any computational efficient procedure under the Gaussian
CCA model considered here under the assumption of planted clique hardness.
The argument has also led to a computational lower bounds for the
sparse PCA problem under the Gaussian spiked covariance model, bridging
the gap between the sparse PCA literature and the computational lower
bounds in \cite{Berthet13} and \cite{wang2014statistical}.

It is of great interest to further investigate if there is some
adaptive and efficient estimator that attains the statistical
optimality established in the current paper under full generality.


\section{Proof of main results}
\label{sec:proof}

This section is devoted to the proof of Theorems~\ref
{thm:upperbound1}--\ref{thm:upperbound2}.
The proof of Theorems \ref{thm:lower-bd-q}--\ref{thm:col-sparse} is
given in the supplementary material \cite{supp2}.




\subsection{Outline of proof and preliminaries}
\label{sec:prelim-proof}
To prove both Theorems \ref{thm:upperbound1} and \ref{thm:upperbound2},
we go through the following three steps:
\begin{longlist}[1.]
\item[1.] We decompose the value of the loss function into multiple terms
which result from different sources;
\item[2.] We derive individual high probability bound for each term in the
decomposition;
\item[3.] We assemble the individual bounds to obtain the desired upper
bounds on the loss and the risk functions.
\end{longlist}

In the rest of this subsection, we carry out these three steps in order.
To facilitate the presentation, we introduce below several important
quantities to be used in the proof.

Recall the effective sparsity $(k_q^u, k_q^v)$ defined in \eqref
{eq:effectivesparsity}.
Let $S_{u}$ be the index set of the rows of $U_{1}$ with the
$k_{q}^{u}$ largest $\ell_2$ norms.
In case $U_{1}$ has no more than $k_q^u$ nonzero rows, we include in
$S_u$ the smallest indices of the zero rows in $U_{1}$ such that $|S_u|
= k_q^u$.
We also define $S_{v}$ analogously.
In what follows, we refer to them as the \emph{effective support sets}.

We define $({U}^*_1, {V}^*_1)$ as a solution to
%
\begin{eqnarray}
\label{eq:UV-oracle} %
 &&\max_{(A,B)}  \Tr
\bigl(A'{\Sigma}_{xy}B\bigr)
\nonumber
\\[-8pt]
\\[-8pt]
\nonumber
&&\qquad
\mbox{s.t.}\quad A'{\Sigma}_xA=B'{
\Sigma}_yB=I_r\mbox{ and } \supp(A) \subset
S_u, \supp(B) \subset S_v.
\end{eqnarray}
In what follows, we refer to them as the \emph{sparse approximations}
to $U_1$ and $V_1$.
By definition, when $q = 0$, ${U}^*_1({V}^*_1)' =U_1 V_1'$, which
can be derived rigorously from Theorem~\ref{thm:sintheta}.

In addition, we define the \emph{oracle estimator} $(\widehat{U}^*_1,
\widehat{V}^*_1)$ as a solution to
%
\begin{eqnarray}
\label{eq:UV-oracle-est} %
 &&\max_{(A,B)}  \Tr
\bigl(A'\wh{\Sigma}_{xy}B\bigr)
\nonumber
\\[-8pt]
\\[-8pt]
\nonumber
&&\qquad\mbox{s.t.} \quad A'\wh{\Sigma}_xA=B'\wh{
\Sigma}_yB=I_r\mbox{ and } \supp(A) = S_u,
\supp(B) = S_v.
\end{eqnarray}
In case the program \eqref{eq:UV-oracle} [or \eqref{eq:UV-oracle-est}]
has multiple global optimizers, we define $({U}^*_1, {V}^*_1)$
[or $(\widehat{U}^*_1, \widehat{V}^*_1)$] by picking an arbitrary one.

\begin{remark}
The introduction of (\ref{eq:UV-oracle}) and (\ref{eq:UV-oracle-est})
is to separate the error brought by not knowing the covariance $\Sigma
_x$ and $\Sigma_y$ and by not knowing the effective supports $S_u$ and
$S_v$. The program (\ref{eq:UV-oracle-est}) assumes known effective
supports but unknown covariance and the program (\ref{eq:UV-oracle})
assumes both known effective supports and known covariance.
\end{remark}

We note that
\[
\bigl({U}^*_1\bigr)_{S_u^c *} = \bigl(\widehat{U}^*_1
\bigr)_{S_u^c *} = 0,\qquad  \bigl({V}^*_1\bigr)_{S_v^c *} =
\bigl(\widehat{V}^*_1\bigr)_{S_v^c *} = 0.
\]
By definition, the matrices $({U}^*_1,{V}^*_1)$ are normalized
with respect to $\Sigma_x$ and $\Sigma_y$,
and $(\widehat{U}^*_1,\widehat{V}^*_1)$ are normalized with respect to
$\widehat{\Sigma}_x$ and $\widehat{\Sigma}_y$. Note the notation
$A_{S*}$ stands for the submatrix of $A$ with rows in $S$ and all columns.

Last but not least, let
%
\begin{equation}
\label{eq:supp-est} \wh{S}_u = \supp(\wh{U}_1),\qquad
\wh{S}_v = \supp(\wh{V}_1).
\end{equation}
By the definition of $(\wh{U}_1, \wh{V}_1)$ in \eqref{eq:pickset}, we
have $|\wh{S}_u| = k_q^u$ and $|\wh{S}_v| = k_q^v$ with probability
one. Remember the minimax rate $\varepsilon_n^2$ defined in (\ref
{eq:DefEpsilon}).

%


\subsection{Loss decomposition} In the first step, we decompose the
loss function into five terms as follows.

\begin{lemma}
\label{lem:lossdecompose}
Assume $\frac{1}{n}
(k_{q}^{u}\log(ep/k_{q}^{u})+k_{q}^{v}\log(em/k_{q}^{v}))<c$ for
sufficiently small $c>0$.
For any constant $C'>0$, there exists a constant $C>0$ only depending
on $M$ and $C'$, such that
%
\begin{eqnarray}
&&\bigl\|\widehat{U}_{1}\widehat{V}_{1}^{\prime}-U_{1}V_{1}^{\prime}
\bigr\| _{\mathrm{F}}^{2}
\nonumber
\\
&&\qquad\leq 3\bigl\| {U}^*_{1}\bigl({V}^*_{1}\bigr)^{\prime
}-U_{1}V_{1}^{\prime}
\bigr\| _{\mathrm{F}}^{2} \label{eq:sparseapproxerror}
\\
&&\qquad\quad{}+3\bigl\| \widehat{U}^*_1\bigl(\widehat{V}^*_1
\bigr)'-{U}^*_{1}%
\bigl({V}^*_{1}
\bigr)^{\prime}\bigr\| _{\mathrm{F}}^{2} \label{eq:oracleloss}
\\
&&\qquad\quad{}-\frac{6C}{\lambda_{r}} \bigl\langle{\Sigma}_{x}{U}_{2}
\Lambda _{2}{V}%
_{2}^{\prime}{
\Sigma}_{y},\widehat{U}^*_1\bigl(\widehat
{V}^*_{1}\bigr)'-\widehat {U}%
_{1}
\widehat{V}_{1}^{\prime} \bigr\rangle\label{eq:bias}
\\
&&\qquad\quad{}+\frac{6C}{\lambda_{r}} \bigl\langle\Sigma_{xy}-\widehat{\Sigma}
_{xy},\widehat{U}^*_1\bigl(
\widehat{V}^*_1\bigr)'-\widehat{U}_{1}
\widehat{V}%
_{1}^{\prime} \bigr\rangle\label{eq:excessloss1}
\\
&&\qquad\quad{}+\frac{6C}{\lambda_{r}} \bigl\langle\widehat{\Sigma}_{x}\widehat
{U}^*_{1}\Lambda_{1}\widehat{V}^*_1{^{\prime}}
\widehat{\Sigma }_{y}-{\Sigma}_{x}{U}%
_{1}
\Lambda_{1}{V}_{1}^{\prime}{\Sigma}_{y},
\widehat{U}^*_1 \bigl(\wh {V}_{1}^{\ast}
\bigr)'-\widehat{U}_{1}\widehat{V}_{1}^{\prime}
\bigr\rangle, \label{eq:excessloss2}
\end{eqnarray}
with probability at least $1-\exp(-C' k_{q}^{u}\log(ep/k_{q}^{u})%
)-\exp(-C' k_{q}^{v}\log(em/k_{q}^{v}))$.
\end{lemma}
\begin{pf}
See Section~\ref{sec:lossdecompose-pf}.
\end{pf}

In particular, Lemma~\ref{lem:lossdecompose} decomposes the total loss
into the sum of the
sparse approximation error in (\ref{eq:sparseapproxerror}), the oracle
loss in (\ref{eq:oracleloss}) which is present even if we have the
oracle knowledge of the effective support sets $S_u$ and $S_v$, the
bias term in (\ref{eq:bias}) caused by the presence
of the residual term $U_2\Lambda_2V_2^{\prime}$ in the CCA structure
(\ref{eq:CCA}) and the
two excess loss terms in (\ref{eq:excessloss1}) and (\ref
{eq:excessloss2}) resulting from the uncertainty about the effective
support sets.
When $q = 0$, the sparse approximation error term \eqref
{eq:sparseapproxerror} vanishes.

\subsection{Bounds for individual terms}
We now state the bounds for the individual terms obtained in Lemma~\ref
{lem:lossdecompose} as five separate lemmas.
The proofs of these lemmas are deferred to Sections \ref
{sec:sparseapproxerror-pf}--\ref{sec:excessloss2-pf}.

\begin{lemma}[(Sparse approximation)]
\label{lem:sparseapproxerror}
Suppose (\ref{eq:ass1}) and (\ref%
{eq:ass2}) hold. There exists a constant $C>0$ only depending on
$M,\kappa,q$, such that
%
\begin{eqnarray}
\bigl\| {U}^*_{1}\bigl({V}^*_{1}\bigr)^{\prime}-U_{1}V_{1}^{\prime
}
\bigr\| _{\mathrm{F}}^{2} &\leq&\frac{Cq}{2-q} \eps_n^2,
%
%
\label{eq:SparseApproClaim1}
\\
\bigl\| {U}^*_{1}\Lambda_{1}\bigl({V}^*_{1}
\bigr)^{\prime
}-U_{1}\Lambda_{1}V_{1}^{\prime}
\bigr\| _{\mathrm{F}}^{2} &\leq& \frac{Cq}{2-q}\lambda^{2}
\eps_n^2. %
\label{eq:SparseApproClaim2}
\end{eqnarray}
\end{lemma}
%

\begin{lemma}[(Oracle loss)]
\label{lem:oracleloss}
Suppose $\frac{1}{n\lambda^{2}} (%
k_{q}^{u}+k_{q}^{v}+\log(ep/k_{q}^{u})+\log(em/k_{q}^{v}) )<c$
and that
(\ref{eq:ass2}) holds for some sufficiently small $c>0$.
For any constant $C'>0$, there exists a constant $C>0$ only depending
on $M,q,\kappa$ and $C'$, such that
%
\begin{equation}
\bigl\| \widehat{U}^*_1\bigl(\widehat{V}^*_1
\bigr)'-{U}^*_{1} \bigl({V}^*_{1}
\bigr)^{\prime}\bigr\| _{\mathrm{F}}^{2}\leq\frac{Cr}{n\lambda
^{2}}%
\biggl[ k_{q}^{u}+k_{q}^{v}+\log
\biggl(\frac{ep}{k_{q}^{u}} \biggr)+\log \biggl(\frac{em}{k_{q}^{v}} \biggr) \biggr],
\label{eq:OracleLossClaim1}
\end{equation}
with probability at least $1-\exp(-C' (k_{q}^{u}+\log
(ep/k_{q}^{u})) )-\exp(-C' (k_{q}^{v}+\break \log
(em/k_{q}^{v})))$.
Moreover, if (\ref{eq:ass1}) also holds,
then with the same probability
%
\begin{equation}
\label{eq:OracleLossClaim2} %
\bigl\| \widehat{U}^*_1
\Lambda_{1}\bigl(\widehat{V}^*_1\bigr)'-{U}^*
_{1}\Lambda_{1}\bigl({V}^*_{1}
\bigr)^{\prime}\bigr\| _{\mathrm{F}}^{2} \leq C 
\lambda^{2} \eps_n^2. 
\end{equation}
\end{lemma}
%

The proof of Lemma~\ref{lem:oracleloss} is given in the supplementary
material \cite{supp2}.
Since $r\leq k_q^u\wedge k_q^v$, (\ref{eq:OracleLossClaim1}) is bounded
above by $C\varepsilon_n^2$. The error bounds in Lemma~\ref
{lem:oracleloss} are due to the estimation error of true covariance
matrices by sample covariance matrices on the subset $S_u\times S_v$.

\begin{lemma}[(Bias)]
\label{lem:bias}
Suppose
$\frac{1}{n}(k_{q}^{u}\log(ep/k_{q}^{u})+k_{q}^{v}\log
(em/k_{q}^{v}))<C_{1}$ for some constant $C_{1}>0$. For any constant $C'>0$,
there exists a constant $C>0$ only depending on $M,\kappa,C_1$ and
$C'$, such that
\begin{eqnarray*}
& &\bigl\llvert \bigl\langle{\Sigma}_{x}{U}_{2}
\Lambda_2{V}_{2}^{\prime
}{\Sigma}%
_{y},
\widehat{U}^*_{1} \bigl(\widehat{V}^*_{1}
\bigr)'-\widehat{U}_{1}\widehat {V}%
_{1}^{\prime}
\bigr\rangle\bigr\rrvert
\\
& &\qquad\leq C\lambda_{r+1} \bigl( \bigl\|\widehat{U}^*_{1} \bigl(
\widehat{V}^*_{1}\bigr)'-{U}_{1}{V}_{1}^{\prime}
\bigr\| _{\mathrm {F}}^2+ \bigl\|U_{1}V^{\prime}_{1}-
\widehat{U}_{1}\widehat {V}_{1}^{\prime}\bigr\|
_{\mathrm{F}}^2 \bigr),
\end{eqnarray*}
with probability at least $1-\exp(-C^{\prime}k_{q}^{u}\log
(ep/k_{q}^{u}))-\exp(-C^{\prime}k_{q}^{v}\log(em/k_{q}^{v}))$.
\end{lemma}
%
The bias in  Lemma~\ref{lem:bias} is $0$ when $U_2\Lambda_2V_2'$ is $0$.

\begin{lemma}[(Excess loss 1)]
\label{lem:excessloss1}
Suppose (\ref{eq:ass1}) holds.
For any constant $C'>0$, there exists a constant $C>0$ only depending
on $M$ and $C'$, such that
\[
\bigl\llvert \bigl\langle\Sigma_{xy}-\widehat{\Sigma}_{xy},
\widehat{U}^*_1 \bigl(\widehat{V}^*_{1}
\bigr)'-\widehat{U}_{1}\widehat {V}_{1}^{\prime
}
\bigr\rangle\bigr\rrvert 
\leq C\lambda\eps_n
\bigl\|\widehat{U}_{1}\widehat{V}_{1}^{\prime}-
\widehat{U}^*_1\bigl(\widehat {V}^*_{1}
\bigr)'\bigr\| _{\mathrm{F}},
\]
with probability at least $1-\exp
(-C'(r(k_{q}^{u}+k_{q}^{v})+k_{q}^{u}\log(ep/k_{q}^{u})+\break k_{q}^{v}\log
(em/ k_{q}^{v})))$.
\end{lemma}
%

\begin{lemma}[(Excess loss 2)]
\label{lem:excessloss2}
Suppose (\ref{eq:ass1}) and (\ref{eq:ass2})
hold.
For any constant \mbox{$C'>0$}, there exists a constant $C>0$ only depending
on $M,\kappa,q$ and $C'$, such that
\begin{eqnarray*}
&&\bigl\llvert \bigl\langle\widehat{\Sigma}_{x}\widehat{U}^*_1
\Lambda_{1} \bigl(\widehat{V}^*_1\bigr)'
\widehat{\Sigma}_{y}-{\Sigma}_{x}{U}_{1}
\Lambda _{1}{V}%
_{1}^{\prime}{
\Sigma}_{y},\widehat{U}^*_1\bigl(\widehat
{V}^*_1\bigr)'-\widehat{U} 
_{1}
\widehat{V}_{1}^{\prime} \bigr\rangle\bigr\rrvert
\\
&&\qquad\leq 
C
\lambda\eps_n\bigl \| \widehat{U}^*_1\bigl(
\widehat{V}^*_1\bigr)'-\widehat{U}_{1}
\widehat{V}%
_{1}^{\prime}\bigr\| _{\mathrm{F}},
\end{eqnarray*}
with probability at least $1-\exp(-C'(k_{q}^{u}+\log
(ep/k_{q}^{u})))-\exp(-C'(k_{q}^{v}+\break \log
(em/k_{q}^{v})))$.
\end{lemma}
%

\subsection{Proof of Theorem \texorpdfstring{\protect\ref{thm:upperbound1}}{1}}
For notational convenience, let
\begin{eqnarray*}
R&=&\bigl\| \widehat{U}_{1}\widehat{V}_{1}^{\prime}-U_{1}V_{1}^{\prime
}
\bigr\|_{\mathrm{F}}, \qquad\theta=\bigl\| {U}^*_{1} \bigl({V}^*_{1}
\bigr)^{\prime}-U_{1}V_{1}^{\prime
}\bigr\|
_{\mathrm{F}},\\
 \delta&=&\bigl\| \widehat{U}^*_1\bigl(\widehat
{V}^*_1\bigr)'-{U}^*_{1}
\bigl({V}^*_{1}\bigr)^{\prime}\bigr\|_{\mathrm{F}}.
\end{eqnarray*}
%
Consider the event such that the conclusions of Lemmas \ref
{lem:lossdecompose}--\ref{lem:excessloss2} hold,
which occurs with probability at least $1-\exp(-C'(k_{q}^{u}+\log
(ep/k_{q}^{u})))-\exp(-C'(k_{q}^{v}+\log(em/k_{q}^{v})))$ according to
the union bound.
On this event,
Lemmas \ref{lem:sparseapproxerror} and \ref{lem:oracleloss}
imply that
\[
\theta^{2}\leq 
C
\varepsilon_n^{2}\quad \mbox{and} \quad\delta^{2} \leq C
\eps_n^2.
\]
%
Moreover,
Lemma~\ref{lem:bias} implies
\[
\biggl\llvert \frac{1}{\lambda_{r}} \bigl\langle{\Sigma }_{x}{U}_{2}
\Lambda _{2}{V}_{2}^{\prime}{\Sigma}_{y},
\widehat{U}^*_1\bigl(\widehat{V}^*_1
\bigr)'-%
\widehat{U}_{1}\widehat{V}_{1}^{\prime}
\bigr\rangle\biggr\rrvert \leq \frac{%
C\lambda_{r+1}}{\lambda}\bigl( R^{2}+
\theta^{2}+\delta^{2} \bigr).
\]
Lemma~\ref{lem:excessloss1} implies
\[
\biggl\llvert \frac{1}{\lambda_{r}} \bigl\langle\Sigma_{xy}-\widehat {
\Sigma}%
_{xy},\widehat{U}^*_1\bigl(
\widehat{V}^*_1\bigr)'-\widehat{U}_{1}
\widehat{V}%
_{1}^{\prime} \bigr\rangle\biggr\rrvert
\leq C\varepsilon_n (R+\theta +\delta),
\]
and Lemma~\ref{lem:excessloss2} implies
\[
\biggl\llvert \frac{1}{\lambda_{r}} \bigl\langle\widehat{\Sigma}%
_{x}
\widehat{U}^*_1\Lambda_{1}\bigl(\widehat{V}^*_1
\bigr){^{\prime}}\widehat {\Sigma }_{y}-{%
\Sigma}_{x}{U}_{1}\Lambda_{1}{V}_{1}^{\prime}{
\Sigma }_{y},U_{1}^{\ast
}\bigl(
\widehat{V}^*_1\bigr)'-\widehat{U}_{1}
\widehat{V}_{1}^{\prime
} \bigr\rangle\biggr\rrvert \leq C
\varepsilon_n (R+\theta+\delta).
\]
Together with Lemma~\ref{lem:lossdecompose}, the above bounds lead to
%
\begin{eqnarray*}
R^{2} &\leq&C\bigl(\theta^{2}+\delta^{2}\bigr)+
\frac{C\lambda_{r+1}}{%
\lambda}\bigl(R^{2}+\theta^{2}+
\delta^{2}\bigr)+C\eps_n (R+\theta+\delta)
\\
&\leq&\frac{C\lambda_{r+1}}{\lambda}R^{2}+C\eps_n R+C
\eps_n ^{2}.
\end{eqnarray*}
Under assumption \eqref{eq:ass2}, we have
$\frac{1}{2}R^{2}\leq C\eps_n R+C\eps_n ^{2}$,
implying
\[
R^{2}\leq C\eps_n ^{2},
\]
for some $C>0$.
We complete the proof by noting that the conditions of Lemmas \ref
{lem:lossdecompose}--\ref{lem:excessloss2} are satisfied under
assumptions \eqref{eq:ass1} and \eqref{eq:ass2}.

\subsection{Proof of Theorem \texorpdfstring{\protect\ref{thm:upperbound2}}{2}}
Recall the definition of $\eps_n$ in \eqref{eq:DefEpsilon}, and let
$C_1$ be the constant in \eqref{eq:TheoremExpAssup1} and \eqref
{eq:TheoremExpAssup2}.
The result of Theorem~\ref%
{thm:upperbound1} implies that we can choose an arbitrarily large
constant $C^{\prime}$
such that $C^{\prime}>C_{1}$. Given~$C'$, there exists a constant~$C$,
by which
we can
bound the risk as follows:
%
\begin{eqnarray}
&&\mathbb{E}_{\Sigma} \bigl\| \widehat{U_{1}V_{1}^{\prime
}}-U_{1}V_{1}^{\prime
}
\bigr\| _{\mathrm{F}}^{2}
\nonumber
\\
&&\qquad\leq \mathbb{E}_\Sigma\bigl[ \bigl\| \widehat{U_{1}V_{1}^{\prime
}}-U_{1}V_{1}^{\prime}
\bigr\|_{\mathrm{F}}^{2} {\mathbf{1}_{\{{\|
\widehat{U_{1}V_{1}^{\prime}}-U_{1}V_{1}^{\prime}\|_{\mathrm
{F}}^{2}\leq C\eps_n ^{2} }\}}} \bigr]
\nonumber
\\
&&\qquad\quad{} +\mathbb{E}_\Sigma\bigl[ \bigl\| \widehat{U_{1}V_{1}^{\prime
}}-U_{1}V_{1}^{\prime}
\bigr\|_{\mathrm{F}}^{2} {\mathbf{1}_{\{{\|
\widehat{U_{1}V_{1}^{\prime}}-U_{1}V_{1}^{\prime}\|_{\mathrm
{F}}^{2} > C\eps_n ^{2} }\}}} \bigr]
\nonumber
\\
&&\qquad\leq C\eps_n ^{2}+\mathbb{E}_\Sigma\bigl[
\bigl( 2\bigl\|\widehat {U_{1}V_{1}^{\prime}}
\bigr\|_{\mathrm{F}}^{2}+2\bigl\| U_{1}V_{1}^{\prime}
\bigr\| _{\mathrm{F}}^{2} \bigr) {\mathbf{1}_{\{{\| \widehat
{U_{1}V_{1}^{\prime}}-U_{1}V_{1}^{\prime}\|_{\mathrm{F}}^{2} > C\eps
_n ^{2} }\}}} \bigr]
\label{eq:TheoremExp1}
\\
&&\qquad\leq C\eps_n ^{2}+ 6M^{2}r
\mathbb{P}_\Sigma\bigl( \bigl\| \widehat{U}_{1}\widehat
{V}_{1}^{\prime}-U_{1}V_{1}^{\prime}
\bigr\|_{\mathrm{F}}^{2}>C\eps_n ^{2} \bigr)
\label{eq:TheoremExp2}
\\
&&\qquad\leq C_{2}\eps_n ^{2}. \label{eq:TheoremExp3}
\end{eqnarray}
Here, inequality (\ref{eq:TheoremExp1}) is due to the triangle
inequality and the fact that
\[
\bigl\{ \bigl\| \widehat{U_{1}V_{1}^{\prime}}-U_{1}V_{1}^{\prime}
\bigr\| _{\mathrm
{F}}^{2}>C\eps_n ^{2} \bigr\}
\subset \bigl\{ \bigl\| \widehat{U}_{1}\widehat{V}_{1}^{\prime}-U_{1}V_{1}^{\prime}
\bigr\| _{\mathrm{F}}^{2}>C\eps_n ^{2} \bigr\}.
\]
In fact, if
$\| \widehat{U}_{1}\widehat{V}_{1}^{\prime}-U_{1}V_{1}^{\prime}\|
_{\mathrm{F}}^{2}\leq C\eps_n ^{2}$, then $\|
\widehat{U}_{1}\widehat{V}_{1}^{\prime}\|_{\mathrm{F}}^{2}\leq
C\eps_n ^{2}+M^{2}r\leq2M^{2}r$. By our definition of the estimator,
this means $\widehat{U_{1}V_{1}^{\prime}}=\widehat{U}_{1}\widehat{V}
_{1}^{\prime}$, which further implies $\| \widehat{%
U_{1}V_{1}^{\prime}}-U_{1}V_{1}^{\prime}\|_{\mathrm{F}}^{2}\leq
C\eps_n ^{2}$.
Inequality (\ref{eq:TheoremExp2}) follows from our definition
of estimator $\widehat{U_{1}V_{1}^{\prime}}$ and (\ref%
{eq:BoundTruth}). The last inequality follows from the conclusion of
Theorem~\ref{thm:upperbound1} and assumptions (\ref%
{eq:TheoremExpAssup1}) and (\ref{eq:TheoremExpAssup2}).
This completes the proof.

\section{Proof of auxiliary results}
\label{sec:aux-proof}

In this section, we prove Lemmas \ref{lem:lossdecompose}--\ref
{lem:sparseapproxerror} and \ref{lem:bias}--\ref{lem:excessloss2} used
in the proof of Theorem~\ref{thm:upperbound1} and \ref{thm:upperbound2}.
The proof of Lemma~\ref{lem:oracleloss} is given in the supplementary
material \cite{supp2}.
Throughout the section, without further notice, $\eps_n^2$ is defined
as in \eqref{eq:DefEpsilon}.

\subsection{A generalized sin-theta theorem and Gaussian quadratic form
with rank constraint} \label{sec:key}
We first introduce two key results used in the proof of Lemmas \ref
{lem:lossdecompose}--\ref{lem:excessloss2}
that might be of independent interest.

The first result is a generalized sin-theta theorem. For the definition
of unitarily invariant norms, we refer the readers to \cite
{Bhatia97,Stewart90}. In particular, both Frobenius norm $\|\cdot\|
_{\mathrm{F}}$ and operator norm $\|\cdot\|_{\mathrm {op}}$ are unitarily invariant.

\begin{theorem}
\label{thm:sintheta} Consider matrices $X,Y\in\mathbb{R}^{p\times m}$.
Let the SVD of $X$ and $Y$ be
\[
X=A_1D_1B_1^{\prime}+A_2D_2B_2^{\prime},\qquad
Y=\widehat{A}_1\widehat {D}_1%
\widehat{B}_1^{\prime}+\widehat{A}_2
\widehat{D}_2\widehat {B}_2^{\prime},
\]
with $D_1=\mathop{\operatorname{diag}}(d_1,\ldots,d_r)$ and
$\widehat{D}_1=
\mathop{\operatorname{diag}}(\widehat{d}_1,\ldots,\widehat{d}_r)$.
Suppose there is a positive constant $\delta\in(0,d_r]$ such that $\|
\widehat{D}_2\|_{\mathrm{op}}\leq d_r-\delta$.
Let $\|\cdot\|$ be any unitarily invariant norm, and
$\varepsilon= \|A_1^{\prime}(X-Y) \| \vee\|(X-Y)B_1 \|$.
Then we have
%
\begin{equation}
\label{eq:MPprincipal}\bigl \|A_1D_1B_1^{\prime}-
\widehat{A}_1\widehat{D}_1\widehat {B}_1^{\prime}
\bigr\| \leq \biggl(\frac{\sqrt{2}(d_1+\widehat{d}_1)}{\delta}+1 \biggr)\varepsilon.
\end{equation}
%
If further there is an absolute constant $\bar{\kappa}\geq1$ such that
$%
d_1\vee\widehat{d}_1\leq\bar{\kappa} d_r$, then there is a constant
$C>0$ only
depending on $\bar{\kappa}$, such that
%
\begin{equation}
\label{eq:MPnew} \bigl\|A_1B_1^{\prime}-
\widehat{A}_1\widehat{B}_1^{\prime}\bigr\| \leq
\frac{C\varepsilon}{\delta}.
\end{equation}
\end{theorem}

\begin{remark}
In addition, when $X$ and $Y$ are positive semi-definite, $A_l=B_l$,
$\wh{A}_l=\wh{B}_l$ for $l=1,2$, we recover the classical Davis--Kahan
sin-theta theorem \cite{Davis70} $\|A_1A_1'-\wh{A}_1\wh{A}_1'\|\leq
C\varepsilon/\delta$ up to a constant multiplier.
\end{remark}

The second result is an empirical process type bound for Gaussian
quadratic forms with rank constraint.

\begin{lemma}
\label{lem:EP} Let $\{Z_{i}\}_{1\leq i\leq n}$ be i.i.d. observations
from $%
N(0,I_{d})$. Then there exist some $C,C^{\prime}>0$, such that for any $t>0$,
\[
\mathbb{P} \Biggl( \sup_{\{K:\| K\| _{\mathrm{F}}\leq1,%
\mathrm{rank}(K)\leq r\}}\Biggl\llvert \Biggl\langle
\frac{1}{n}%
\sum_{i=1}^{n}Z_{i}Z_{i}^{\prime}-I_{d},K
\Biggr\rangle\Biggr\rrvert >t \Biggr) \leq\exp\bigl(C^{\prime}rd-Cn
\bigl(t^{2}\wedge t\bigr)\bigr).
\]
\end{lemma}

The proofs of Theorem~\ref{thm:sintheta} and Lemma~\ref{lem:EP} are
given in the supplementary material~\cite{supp2}.

\subsection{Proof of Lemma \texorpdfstring{\protect\ref{lem:lossdecompose}}{1}}
\label{sec:lossdecompose-pf}

Recall the definition of $(S_u, S_v)$ and $(\wh{S}_u, \wh{S}_v)$ in
Section~\ref{sec:prelim-proof}.
From here on, let
%
\begin{equation}
\label{eq:Tu-Tv} T_{u}=S_{u}\cup\widehat{S}_{u}\quad
\mbox{and}\quad T_{v}=S_{v}\cup\widehat{S}_{v}.
\end{equation}
The proof of Lemma~\ref{lem:lossdecompose} depends on the following two
technical results. Their proofs are given in the supplementary material
\cite{supp2}.

\begin{lemma}
\label{lem:linearloss} For matrices $A,B,E,F$ and a diagonal matrix $%
D=(d_{l})_{1\leq l\leq r}$ with $d_{1}\geq d_{2}\geq\cdots\geq d_{r}>0$
and $%
A^{\prime}A=B^{\prime}B=E^{\prime}E=F^{\prime}F=I_r$, we have
\[
\frac{d_{r}}{2}\bigl\| AB^{\prime}-EF^{\prime}\bigr\| _{\mathrm{F}%
}^{2}
\leq \bigl\langle ADB^{\prime},AB^{\prime}-EF^{\prime} \bigr
\rangle \leq\frac{d_{1}}{2}\bigl\| AB^{\prime}-EF^{\prime}\bigr\|
_{%
\mathrm{F}}^{2}.
\]
\end{lemma}

\begin{lemma}
\label{lem:diff}
Under the assumption of Lemma~\ref{lem:lossdecompose},
for any constant $C'>0$, there exists a constant $C>0$ only depending
on $M$ and $C'$, such that for any matrix $A$ supported on the
$T_{u}\times T_{v}$,
we have
\[
{C}^{-1}
\|A\|_{\mathrm{F}}^2\leq \bigl\|\widehat{\Sigma}_{x}^{1/2}
A \widehat{\Sigma}_{y}^{1/2}\bigr\|_{\mathrm{F}}^{2}
\leq
C\| A\| _{\mathrm{F}}^{2},
\]
with probability at least $1-\exp(-C^{\prime}k_{q}^{u}\log
(ep/k_{q}^{u}))-\exp(-C^{\prime}k_{q}^{v}\log(em/k_{q}^{v}))$.
\end{lemma}

\begin{pf*}{Proof of Lemma~\ref{lem:lossdecompose}}
First of all, the triangle inequality and Jensen's inequality together
lead to
%
\begin{eqnarray}
\label{eq:losssurrogate} %
&& \bigl\| \widehat{U}_{1}
\widehat{V}_{1}^{\prime}-U_{1}V_{1}^{\prime
}
\bigr\| _{\mathrm{F}}^{2}
\nonumber
\\[-8pt]
\\[-8pt]
\nonumber
&&\qquad \leq3\bigl( \bigl\| \widehat{U}_{1}^{\ast}\widehat{V}_{1}^{\ast\prime
}-U_{1}^{\ast}V_{1}^{\ast\prime}
\bigr\| _{\mathrm{F}}^{2} + \bigl\| \widehat {U}_{1}
\widehat{V}_{1}^{\prime}-\widehat{U}_{1}^{\ast }
\widehat {V}_{1}^{\ast\prime}\bigr\| _{\mathrm{F}}^{2} +\bigl \|
U_{1}^{\ast
}V_{1}^{\ast\prime}-U_{1}V_{1}^{\prime }
\bigr\| _{\mathrm{F}}^{2} \bigr).
\end{eqnarray}
Now, it remains to bound $\| \widehat{U}_{1}\widehat{V}_{1}^{\prime
}-\widehat{U}_{1}^{\ast}\widehat{V}_{1}^{\ast\prime}\| _{\mathrm{%
F}}^{2}$.
To this end, we have
%
\begin{eqnarray}
&& \bigl\| \widehat{U}_{1}^{\ast}\widehat{V}_{1}^{\ast\prime}-
\widehat{U 
}_{1}\widehat{V}_{1}^{\prime}
\bigr\| _{\mathrm{F}}^{2}
\nonumber
\\
&&\qquad\leq C\bigl\| \widehat{\Sigma}_{x}^{1/2}\bigl(
\widehat{U}_{1}^{\ast}%
\widehat{V}_{1}^{\ast\prime}-
\widehat{U}_{1}\widehat {V}_{1}^{\prime}\bigr)
\widehat{\Sigma}_{y}^{1/2}\bigr\| _{\mathrm{F}}^{2}
\label{eq:sigma1}
\\
&&\qquad\leq\frac{2C}{\lambda_{r}} \bigl\langle\widehat{\Sigma}_{x}^{1/2}
\widehat{U}_{1}^{\ast}\Lambda_{1}
\widehat{V}_{1}^{\ast\prime
}\widehat{%
\Sigma}_{y}^{1/2},\widehat{\Sigma}_{x}^{1/2}
\bigl(\widehat{U}_{1}^{\ast} 
\widehat{V}_{1}^{\ast\prime}-
\widehat{U}_{1}\widehat {V}_{1}^{\prime}\bigr)
\widehat{\Sigma}_{y}^{1/2} \bigr
\rangle\label{eq:linearloss}
\\
&&\qquad = \frac{2C}{\lambda_{r}} \bigl\langle\widehat{\Sigma }_{x}\widehat{U}
_{1}^{\ast}\Lambda_{1}
\widehat{V}_{1}^{\ast\prime}\widehat{\Sigma }_{y},%
\widehat{U}_{1}^{\ast}\widehat{V}_{1}^{\ast\prime}-
\widehat{U}_{1} 
\widehat{V}_{1}^{\prime}
\bigr\rangle
\nonumber
\\
&&\qquad = \frac{2C}{\lambda_{r}} \bigl\langle\widehat{\Sigma }_{x}\widehat{U}
_{1}^{\ast}\Lambda_{1}
\widehat{V}_{1}^{\ast\prime}\widehat{\Sigma }_{y}-%
\widehat{\Sigma}_{xy},\widehat{U}_{1}^{\ast}
\widehat{V}_{1}^{\ast
\prime}-%
\widehat{U}_{1}
\widehat{V}_{1}^{\prime} \bigr\rangle
\nonumber
\\
&&\qquad\quad{} +\frac{2C}{\lambda_{r}}%
 \bigl\langle\widehat{\Sigma}_{xy},
\widehat{U}_{1}^{\ast}\widehat {V}%
_{1}^{\ast\prime}-
\widehat{U}_{1}\widehat{V}_{1}^{\prime} \bigr\rangle
\nonumber
\\
&&\qquad\leq\frac{2C}{\lambda_{r}} \bigl\langle\widehat{\Sigma }_{x}\widehat
{U}%
_{1}^{\ast}\Lambda_{1}
\widehat{V}_{1}^{\ast\prime}\widehat{\Sigma }_{y}-%
\widehat{\Sigma}_{xy},\widehat{U}_{1}^{\ast}
\widehat{V}_{1}^{\ast
\prime}-%
\widehat{U}_{1}
\widehat{V}_{1}^{\prime} \bigr\rangle\label{eq:aggregdef}
\\
&&\qquad = \frac{2C}{\lambda_{r}} \bigl\langle\widehat{\Sigma }_{x}\widehat{U}
_{1}^{\ast}\Lambda_{1}
\widehat{V}_{1}^{\ast\prime}\widehat{\Sigma }_{y}-{%
\Sigma}_{x}{U}_{1}\Lambda_{1}{V}_{1}^{\prime}{
\Sigma}_{y},\widehat {U}%
_{1}^{\ast}
\widehat{V}_{1}^{\ast^{\prime}}-\widehat{U}_{1}\widehat{V}
_{1}^{\prime} \bigr\rangle\label{eq:biasstructure1}
\\
&&\qquad\quad{} - \frac{2C}{\lambda_{r}} \bigl\langle{\Sigma }_{x}{U}_{2}
\Lambda _{2}{V}%
_{2}^{\prime}{
\Sigma}_{y},\widehat{U}_{1}^{\ast}\widehat
{V}_{1}^{\ast
^{\prime}}-\widehat{U}_{1}
\widehat{V}_{1}^{\prime} \bigr\rangle
\nonumber
\\
&&\qquad\quad{} + \frac{2C}{\lambda_{r}} \bigl\langle\Sigma_{xy}-\widehat {
\Sigma}_{xy},%
\widehat{U}_{1}^{\ast}
\widehat{V}_{1}^{\ast^{\prime}}-\widehat{U}_{1} 
\widehat{V}_{1}^{\prime} \bigr\rangle.
\nonumber
\end{eqnarray}
%
Here, (\ref{eq:sigma1}) is implied by Lemma~\ref{lem:diff} and (\ref
{eq:linearloss}) is implied by Lemma~\ref{lem:linearloss}.
To see (\ref{eq:aggregdef}), we note $(\wh{U}_1, \wh{V}_1)$ is the
solution to \eqref{eq:pickset}, and so
$\mathop{\sf Tr}(\widehat{U}_{1}^{\prime}\widehat{\Sigma
}_{xy}\widehat{V}_1)
\geq
\mathop{\sf Tr}((\widehat{U}_{1}^{\ast})'\widehat{\Sigma}_{xy} \wh
{V}_{1}^{\ast})$,
or equivalently
\[
\bigl\langle\widehat{\Sigma}_{xy},\widehat{U}_{1}^{\ast}
\widehat {V}%
_{1}^{\ast\prime}-\widehat{U}_{1}
\widehat{V}_{1}^{\prime} \bigr\rangle \leq0.
\]
Equality (\ref{eq:biasstructure1}) comes from
the CCA structure (\ref{eq:CCA}) and (\ref{eq:splitCCA}). Combining
(\ref{eq:losssurrogate})--(\ref%
{eq:biasstructure1}) and rearranging the terms, we obtain the desired
result. 
\end{pf*}

\subsection{Proof of Lemma \texorpdfstring{\protect\ref{lem:sparseapproxerror}}{2}}
\label{sec:sparseapproxerror-pf}

The major difficulty in proving the lemma lies in the presence of the
residual structure $U_2 \Lambda_2 V_2'$ in \eqref{eq:splitCCA} and the
possible nondiagonality of covariance matrices $\Sigma_x$ and $\Sigma_y$.
To overcome the difficulty, we introduce intermediate matrices
$(\wt{U}_{1},\wt{V}_{1})$ defined as follows.
First, we write the SVD of $(\Sigma
_{xS_{u}S_{u}})^{1/2}{U}_{1S_{u}\ast}\Lambda_{1}({V}_{1S_{v}\ast
})^{\prime}(\Sigma_{yS_{v}S_{v}})^{1/2}$ as
%
\begin{equation}
(\Sigma_{xS_{u}S_{u}})^{1/2}{U}_{1S_{u}\ast}\Lambda_{1}({V}%
_{1S_{v}\ast})^{\prime}(
\Sigma_{yS_{v}S_{v}})^{1/2}=P\wt{%
\Lambda}_{1}Q^{\prime},
\label{eq:intermSVD}
\end{equation}
and let
$\wt{U}_{1}^{S_{u}} = ( \Sigma_{xS_{u}S_{u}} )^{-1/2} P$ and
$\wt{V}_{1}^{S_{v}} = ( \Sigma_{yS_{v}S_{v}} )^{-1/2} Q$.
Finally, we define $\wt{U}_1\in\reals^{p\times r}$ and $\wt{V}_1
\in
\reals^{m\times r}$ by setting
%
\begin{equation}
\label{eq:UV-intermediate} \qquad(\wt{U}_1)_{S_u *} = \wt{U}_{1}^{S_{u}},\qquad
(\wt{U}_1)_{S_u^c *} = 0,\qquad (\wt{V}_1)_{S_v *}
= \wt{V}_{1}^{S_{v}},\qquad (\wt{V}_1)_{S_v^c *}
= 0.
\end{equation}
By definition, we have ${U}_{1S_{u}\ast}\Lambda_{1}({V}_{1S_{u}\ast
})^{\prime}=\wt{U}_{1S_{u}\ast}\wt{\Lambda}({\wt{V}}_{1S_{u}\ast}
)^{\prime}$.
Last but not least, we define
%
\begin{eqnarray}\quad\hspace*{4pt}
\label{eq:Xidef} \Xi&=& P\wt{\Lambda}_{1}Q^{\prime}
\nonumber
\\[-8pt]
\\[-8pt]
\nonumber
&&{} +
\bigl(I-PP^{\prime}\bigr) (\Sigma_{xS_{u}S_{u}}%
)^{-1/2}
\Sigma_{xS_{u}\ast}U_{2}\Lambda_{2}V_{2}^{\prime}
\Sigma _{y\ast S_{v}}(\Sigma_{yS_{v}S_{v}})^{-1/2}
\bigl(I-QQ^{\prime}\bigr).
\end{eqnarray}
We now summarize the key properties of the $\wt{U}_1, \wt{V}_1$ and
$\wt
{\Lambda}_1$ matrices in the following two lemmas, the proofs of which
are given in the supplementary material~\cite{supp2}.

\begin{lemma}
\label{lem:CCA2}
Let $P, Q$ and $\Xi$ be defined in \eqref{eq:intermSVD} and \eqref{eq:Xidef}.
Then we have:
\begin{longlist}[1.]
\item[1.] The column vectors of $P$ and $Q$ are the $r$ leading left and
right singular vectors of $\Xi$.
\item[2.] The first and the $r$th singular values $\wt{\lambda}_1$ and
$\wt
{\lambda}_r$ of $\Xi$ satisfy $1.1\kappa\lambda\geq\wt{\lambda}_1
\geq\wt{\lambda}_{r}\geq0.9\lambda$, and the $(r+1)$th singular value
$\wt{\lambda}_{r+1}\leq c\lambda$ for some sufficiently small
constant $c>0$.
\item[3.] The column vectors of $\Sigma_{x}^{1/2}\wt{U}_{1}$ and $\Sigma
_{y}^{1/2}\wt{V}_{1}
$ are the $r$ leading left and right singular vectors of $\Sigma
_{x}^{1/2}\wt{U}_{1}\wt{%
\Lambda}_{1}\wt{V}_{1}^{\prime}\Sigma_{y}^{1/2}$.
\end{longlist}
\end{lemma}

\begin{lemma}
\label{lem:lq3}
For some constant $C > 0$,
\[
\bigl\| \wt{U}_1' \Sigma_x U_{2}\bigr\|
_{\mathrm{F}}^{2} \leq C\| U_{1S_{u}^{c}\ast}\| _{\mathrm{F}}^{2}\quad
\mbox{and} \quad\bigl\| \wt{V}_1' \Sigma_y
V_{2}\bigr\| _{\mathrm{F}}^{2} \leq C\| V_{1S_{v}^{c}\ast}\|
_{\mathrm{F}}^{2}.
\]
%
\end{lemma}

In what follows, we prove claims \eqref{eq:SparseApproClaim1} and
\eqref
{eq:SparseApproClaim2} in order.

\begin{pf*}{Proof of (\ref{eq:SparseApproClaim1})}
By the  triangle inequality,
%
\begin{equation}
\bigl\| U_{1}^{\ast}V_{1}^{\ast\prime}-U_{1}V_{1}^{\prime}
\bigr\| _{\mathrm{F}}\leq\bigl\| U_{1}^{\ast}V_{1}^{\ast\prime}-
\wt{U}_{1}%
\wt{V}_{1}^{\prime}\bigr\|
_{\mathrm{F}}+\bigl\| \wt{U}_{1}\wt{V}%
_{1}^{\prime}-U_{1}V_{1}^{\prime}
\bigr\| _{\mathrm{F}}. \label{eq:piece1}
\end{equation}
It is sufficient to bound each of the two terms on the right-hand side.

\textit{$1^\circ$ Bound for $\| \wt{U}_{1}\wt{V}_{1}^{\prime
}-U_{1}V_{1}^{\prime}\| _{\mathrm{F}}$}.
Since the smallest eigenvalues of $\Sigma_x$ and $\Sigma_y$ are bounded
from below by some absolute positive constant,
\[
\bigl\| \wt{U}_{1}\wt{V}_{1}^{\prime}-U_{1}V_{1}^{\prime}
\bigr\|_{\mathrm{F}} \leq C\bigl\| \Sigma_{x}^{1/2} \bigl(
\wt{U}_{1}\wt{V}_{1}^{\prime
}-U_{1}V_{1}^{\prime}
\bigr)\Sigma_{y}^{1/2}\bigr\|_{\mathrm{F}}.
\]
By Lemma~\ref{lem:CCA2}, $\Sigma_{x}^{1/2}\wt{U}_{1}$ and $\Sigma
_{y}^{1/2}\wt{V}_{1}$ collect the $r$ leading left and right singular
vectors of $\Sigma_{x}^{1/2}\wt{U}_{1}\wt{\Lambda}_{1}\wt
{V}_{1}^{\prime}\Sigma_{y}^{1/2}$, and by (\ref{eq:CCA}), $\Sigma
_{x}^{1/2}{U}_{1}$ and $\Sigma_{y}^{1/2}{V}_{1} $ collect the $r$
leading left and right singular vectors of $\Sigma_{x}^{1/2}{U}_{1}{%
\Lambda}_{1}{V}_{1}^{\prime}\Sigma_{y}^{1/2}$.
Thus, Theorem~\ref{thm:sintheta} implies
\[
\bigl\| \Sigma_{x}^{1/2}\bigl(\wt{U}_{1}
\wt{V}_{1}^{\prime
}-U_{1}V_{1}^{\prime}
\bigr)\Sigma_{y}^{1/2}\bigr\| _{\mathrm{F}}\leq
\frac
{C}{\lambda} \bigl\| \Sigma_{x}^{1/2}\bigl(
\wt{U}_{1}\wt{\Lambda}_{1}\wt{V}_{1}^{\prime
}-U_{1}
\Lambda_{1}V_{1}^{\prime}\bigr)\Sigma_{y}^{1/2}
\bigr\| _{\mathrm{F}}.
\]
The right-hand side of the above inequality is bounded as
\begin{eqnarray}
&& \label{eq:zhaoren} \bigl\| \wt{U}_{1}\wt{\Lambda}_{1}
\wt{V}_{1}^{\prime
}-U_{1}\Lambda_{1}V_{1}^{\prime}
\bigr\| _{\mathrm{F}} 
\nonumber\\
&& \qquad\leq\bigl\| \wt{U}_{1S_{u}\ast}\wt{\Lambda}_{1}(
\wt{V}_{1S_{v}\ast
})^{\prime}-U_{1S_{u}\ast}\Lambda
_{1}(V_{1S_{v}\ast})^{\prime}\bigr\|_{\mathrm{F}}
+
\bigl\|U_{1S_{u}^{c}\ast}\Lambda_{1}(V_{1S_{v}\ast})^{\prime
}\bigr\|
_{\mathrm{F}}
\nonumber
\\[-8pt]
\\[-8pt]
\nonumber
&&\qquad\quad{} +\bigl \| U_{1S_{u}\ast}\Lambda_{1}(V_{1S_{v}^{c}\ast})^{\prime}
\bigr\| _{\mathrm{F}} +\bigl \|U_{1S_{u}^{c}\ast}
\Lambda_{1}(V_{1S_{v}^{c}\ast}
)^{\prime}\bigr\| _{\mathrm{F}}
\nonumber
\\
\nonumber
&&\qquad\leq C\lambda\bigl(\| U_{1S_{u}^{c}\ast}\|_{\mathrm{F}}+\|
V_{1S_{v}^{c}\ast}\| _{\mathrm{F}}\bigr).
\nonumber
\end{eqnarray}
Here, the last inequality is due to \eqref{eq:intermSVD} and \eqref
{eq:UV-intermediate}.
For the last term, a similar argument to that used in Lemma 7 of \cite
{CMW13a} leads to
%
\begin{eqnarray}
\label{eq:lqsparse} %
\bigl \| U_{1S_{u}^{c}\ast}
\bigr\|_{\mathrm{F}}^2 & \leq&\frac
{Cq}{2-q}k_{q}^{u}
\bigl(s_{u}/k_{q}^{u}\bigr)^{2/q} \leq
\frac{Cq}{2-q} \eps_n^2,
\nonumber
\\[-8pt]
\\[-8pt]
\nonumber
\|V_{1S_{v}^{c}\ast}\| _{\mathrm{F}}^2 & \leq& \frac{Cq}{2-q}
k_{q}^{v}\bigl(s_{v}/k_{q}^{v}
\bigr)^{2/q} \leq\frac{Cq}{2-q} \eps_n^2,
\end{eqnarray}
where the last inequalities in both displays are due to \eqref
{eq:x-q-u}--\eqref{eq:effectivesparsity}.
Therefore, we obtain
%
\begin{equation}
\bigl\| \wt{U}_{1}\wt{V}_{1}^{\prime}-U_{1}V_{1}^{\prime}
\bigr\| _{\mathrm{F}}^{2}\leq\frac{Cq}{2-q} \eps_n^2.
\label{eq:piece2}
\end{equation}

\textit{$2^\circ$ Bound for $\| U_{1}^{\ast}V_{1}^{\ast
\prime}-\wt{U}_{1}\wt{V}_{1}^{\prime}\| _{\mathrm{F}}$}.
Let $\lambda_{r+1}^{\ast}$ denote the $(r+1)$th singular value of
$(\Sigma_{xS_{u}S_{u}})
^{-1/2}\Sigma_{xyS_{u}S_{v}}(\Sigma_{yS_{v}S_{v}})^{-1/2}$.
Then we have
%
\begin{eqnarray}
\label{eq:bound2ndSparseAppr} &&\bigl\| U_{1}^{\ast}V_{1}^{\ast\prime}-
\wt{U}_{1}\wt{V}_{1}^{\prime}\bigr\| _{\mathrm{F}}
\nonumber
\\
&&\qquad = \bigl\| U_{1S_{u}\ast}^{\ast}\bigl(V_{1S_{v}\ast}^{\ast}
\bigr)^{\prime}-\wt {U}_{1S_{u}\ast}(\wt{V}_{1S_{v}\ast})^{\prime}
\bigr\| _{\mathrm{F}}
\nonumber
\\[-8pt]
\\[-8pt]
\nonumber
&&\qquad\leq C\bigl\| (\Sigma_{xS_{u}S_{u}})^{1/2}\bigl[U_{1S_{u}\ast
}^{\ast}
\bigl(V_{1S_{v}\ast}^{\ast}\bigr)^{\prime}-\wt{U}_{1S_{u}\ast}(
\wt {V}_{1S_{v}\ast})^{\prime}\bigr] (\Sigma_{yS_{v}S_{v}})^{1/2}
\bigr\|_{\mathrm{F}}
\\
&&\qquad\leq\frac{C
\| (\Sigma_{xS_{u}S_{u}})^{-1/2}\Sigma_{xyS_{u}S_{v}}(\Sigma
_{yS_{v}S_{v}})^{-1/2}-\Xi\| _{\mathrm{F}} }{
\wt{\lambda}_{r}-\lambda_{r+1}^{\ast}}.\nonumber
\end{eqnarray}
Here, the first equality holds since both $U_{1}^{\ast}V_{1}^{\ast
\prime}$ and $\wt{U}_{1}\wt{V}_{1}^{\prime}$ are supported on the
$S_u\times S_v$ submatrix.
Noting that by the discussion before \eqref{eq:UV-oracle}, \eqref
{eq:UV-intermediate} and Lemma~\ref{lem:CCA2}, $((\Sigma
_{xS_{u}S_{u}})^{1/2}U_{1S_{u}\ast}^{\ast}, (\Sigma
_{yS_{v}S_{v}})^{1/2}V_{1S_{v}\ast}^{\ast})$ and $((\Sigma
_{xS_{u}S_{u}})^{1/2}\wt{U}_{1S_{u}\ast}, \break (\Sigma
_{yS_{v}S_{v}})^{1/2}\wt{V}_{1S_{v}\ast})$ collect the leading left
and right singular vectors of\break $(\Sigma_{xS_{u}S_{u}})^{-1/2}\Sigma
_{xyS_{u}S_{v}} (%
\Sigma_{yS_{v}S_{v}})^{-1/2}$ and $\Xi$, respectively,
we obtain the last inequality by applying (\ref{eq:MPnew}) in
Theorem~\ref{thm:sintheta}.
In what follows, we derive upper bound for the numerator and lower
bound for the denominator in \eqref{eq:bound2ndSparseAppr} in order.

\textit{Upper bound for $\| (\Sigma_{xS_{u}S_{u}})^{-1/2}\Sigma
_{xyS_{u}S_{v}}(\Sigma_{yS_{v}S_{v}})^{-1/2}-\Xi\| _{\mathrm{F}}$}.
First, we decompose $\Sigma_{xyS_{u}S_{v}}$ as
%
\begin{eqnarray}\label{eq:decompSigma_xy}
\Sigma_{xyS_{u}S_{v}} & =& \Sigma_{xS_{u}\ast}\bigl(U_{1}\Lambda
_{1}V_{1}^{\prime}+U_{2}
\Lambda_{2}V_{2}^{\prime}\bigr)\Sigma_{y\ast
S_{v}}
\nonumber
\\
& = &\Sigma_{xS_{u}S_{u}}U_{1S_{u}\ast}\Lambda_{1}V_{1S_{v}\ast
}^{\prime
}
\Sigma_{yS_{v}S_{v}}+\Sigma_{xS_{u}S_{u}}U_{1S_{u}\ast}\Lambda
_{1}V_{1S_{v}^{c}\ast}^{\prime}\Sigma_{yS_{v}^{c}S_{v}}
\\
&&{} +\Sigma_{xS_{u}S_{u}^{c}}U_{1S_{u}^{c}\ast}\Lambda_{1}V_{1
}^{\prime}
\Sigma_{y\ast S_{v}}+\Sigma_{xS_{u}\ast}U_{2}\Lambda
_{2}V_{2}^{\prime}\Sigma_{y\ast S_{v}}.
\nonumber
\end{eqnarray}
Then \eqref{eq:decompSigma_xy}, \eqref{eq:Xidef} and \eqref
{eq:intermSVD} jointly imply that
\begin{eqnarray*}
&& \bigl\| (\Sigma_{xS_{u}S_{u}})^{-1/2}\Sigma_{xyS_{u}S_{v}}(\Sigma
_{yS_{v}S_{v}})^{-1/2}-\Xi\bigr\| _{\mathrm{F}}
\\
&&\qquad \leq\bigl\| (\Sigma_{xS_{u}S_{u}}) ^{-1/2}\Sigma_{xS_{u}S_{u}^{c}}U_{1S_{u}^{c}\ast}
\Lambda_{1}V_{1
}^{\prime}\Sigma_{y\ast S_{v}}(
\Sigma_{yS_{v}S_{v}}) ^{-1/2}\bigr\| _{\mathrm{F}}
\\
& &\qquad\quad{}+\bigl\| (\Sigma_{xS_{u}S_{u}})^{1/2}U_{1S_{u}\ast}\Lambda
_{1}V_{1S_{v}^{c}\ast}^{\prime}\Sigma_{yS_{v}^{c}S_{v}}(\Sigma
_{yS_{v}S_{v}})^{-1/2}\bigr\| _{\mathrm{F}}
\nonumber
\\
&&\qquad\quad{} +\bigl\| PP^{\prime}(\Sigma_{xS_{u}S_{u}})^{-1/2}\Sigma
_{xS_{u}\ast}U_{2}\Lambda_{2}V_{2}^{\prime}
\Sigma_{y\ast S_{v}}(%
\Sigma_{yS_{v}S_{v}})^{-1/2}
\bigl(I-QQ^{\prime}\bigr)\bigr\| _{\mathrm{F}}
\\
&&\qquad\quad{} + \bigl\| (\Sigma_{xS_{u}S_{u}})^{-1/2}\Sigma _{xS_{u}\ast}U_{2}
\Lambda_{2}V_{2}^{\prime}\Sigma_{y\ast S_{v}}(%
\Sigma_{yS_{v}S_{v}})^{-1/2}QQ^{\prime}\bigr\| _{\mathrm{F}%
}
\\
&&\qquad \leq C \lambda\bigl(\| U_{1S_{u}^{c}\ast}\| _{\mathrm{F}}+\| V_{1S_{v}^{c}\ast}\|
_{\mathrm{F}}\bigr)
\\
& &\qquad\quad{}+ C\lambda_{r+1} \bigl(\bigl\| P^{\prime}(\Sigma
_{xS_{u}S_{u}})^{-1/2}\Sigma _{xS_{u}\ast}U_{2}\bigr\|
_{\mathrm{F}}+\bigl\| Q^{\prime}(\Sigma _{yS_{v}S_{v}})^{-1/2}
\Sigma_{yS_{v}\ast}V_{2}\bigr\|_{\mathrm{F}}\bigr)
\\
&&\qquad = C \lambda\bigl(\| U_{1S_{u}^{c}\ast}\| _{\mathrm{F}}+\| V_{1S_{v}^{c}\ast}\|
_{\mathrm{F}}\bigr) + C\lambda_{r+1}\bigl(\bigl\| \wt{U}_1'
\Sigma_x U_{2}\bigr\| _{\mathrm{F}} + \bigl\| \wt
{V}_1' \Sigma_y V_{2}\bigr\|
_{\mathrm{F}}\bigr).
\end{eqnarray*}
Here, the last equality is due to the definition \eqref{eq:UV-intermediate}.
The last display, together with~\eqref{eq:lqsparse} and Lemma~\ref
{lem:lq3}, leads to
%
\begin{equation}
\bigl\| (\Sigma_{xS_{u}S_{u}})^{-1/2}\Sigma_{xyS_{u}S_{v}}(%
\Sigma_{yS_{v}S_{v}})^{-1/2}-\Xi\bigr\| _{\mathrm{F}}^{2}\leq
\frac{Cq}{2-q} \lambda^{2}\eps_n^2.
\label{eq:usedlater4}
\end{equation}

\textit{Lower bound for $\wt{\lambda}_{r}-\lambda_{r+1}^{\ast}$}.
The\vspace*{1pt} bound \eqref{eq:usedlater4}, together with Weyl's inequality (\cite
{golub1996matrix},
page 449 and Hoffman--Wielant inequality \cite{tao12},
page 63) implies
%
\begin{eqnarray}
\label{eq:boundsLamdaTilda}\qquad %
&&\bigl |\lambda_{r+1}^{\ast}-\wt{
\lambda}_{r+1}\bigr| \vee\bigl\| \Lambda_{1}^{\ast}-\wt{
\Lambda}_{1}\bigr\|_{\mathrm{F}}
\nonumber
\\[-8pt]
\\[-8pt]
\nonumber
&&\qquad \leq\bigl\| (\Sigma_{xS_{u}S_{u}})^{-1/2}\Sigma_{xyS_{u}S_{v}}(%
\Sigma_{yS_{v}S_{v}})^{-1/2}-\Xi\bigr\| _{\mathrm{F}} 
\leq C \sqrt{\frac{q}{2-q}} \lambda\eps_n %
\leq0.1\lambda. %
\end{eqnarray}
Together with Lemma~\ref{lem:CCA2}, it further implies
%
\begin{equation}
\label{eq:lambda-lower} \wt{\lambda}_{r}-\lambda_{r+1}^{\ast}
\geq\wt{\lambda}_{r}-\wt {\lambda }%
_{r+1}-\bigl|\wt{
\lambda}_{r+1}-\lambda_{r+1}^{\ast}\bigr|\geq0.7\lambda.
\end{equation}
Combining \eqref{eq:bound2ndSparseAppr}, \eqref{eq:usedlater4} and
\eqref{eq:lambda-lower}, we obtain
%
\begin{equation}
\bigl\| \wt{U}_{1}\wt{V}_{1}^{\prime}-U_{1}^{\ast}V_{1}^{\ast\prime
}
\bigr\| _{\mathrm{F}}^{2}\leq\frac{Cq}{2-q} \eps_n^2.
\label{eq:piece3}
\end{equation}
The proof of \eqref{eq:SparseApproClaim1} is completed by combining
\eqref{eq:piece1}, \eqref{eq:piece2} and \eqref{eq:piece3}.
\end{pf*}

\begin{pf*}{Proof of (\ref{eq:SparseApproClaim2})}
Note that
\begin{eqnarray*}
&&\bigl\| U_{1}^{\ast}\Lambda_{1}V_{1}^{\ast\prime}-{U}_{1}
\Lambda _{1}{V}_{1}^{\prime}\bigr\| _{\mathrm{F}}
\\
&&\qquad\leq\bigl\| U_{1}^{\ast}\Lambda_{1}V_{1}^{\ast\prime}-
\wt{U}_{1}%
\wt{\Lambda}_{1}\wt{V}_{1}^{\prime}
\bigr\| _{\mathrm{F}}+\bigl\| \wt{U}_{1}\wt{\Lambda}_{1}
\wt{V}_{1}^{\prime}-{U}_{1}\Lambda
_{1}{V}%
_{1}^{\prime}\bigr\| _{\mathrm{F}}
\\
&&\qquad\leq\bigl\| U_{1}^{\ast}\Lambda_{1}^{\ast}V_{1}^{\ast\prime}-
\wt{%
U}_{1}\wt{\Lambda}_{1}\wt{V}_{1}^{\prime}
\bigr\| _{\mathrm{F}%
}+\bigl\| \wt{U}_{1}\wt{\Lambda}_{1}
\wt{V}_{1}^{\prime}-{U}%
_{1}
\Lambda_{1}{V}_{1}^{\prime}\bigr\| _{\mathrm{F}}
\\
&&\qquad\quad{} +C\bigl\| \Lambda_{1}^{\ast}-\wt{\Lambda}_{1}\bigr\|
_{\mathrm{F}}+C\| \wt{\Lambda}_{1}-\Lambda_{1}\|
_{\mathrm{F}}
\\
&&\qquad\leq\bigl\| U_{1}^{\ast}\Lambda_{1}^{\ast}V_{1}^{\ast\prime}-
\wt{%
U}_{1}\wt{\Lambda}_{1}\wt{V}_{1}^{\prime}
\bigr\| _{\mathrm{F}%
}+C^{\prime}\bigl\| \wt{U}_{1}\wt{
\Lambda}_{1}\wt{V}_{1}^{\prime}-{U}%
_{1}
\Lambda_{1}{V}_{1}^{\prime}\bigr\| _{\mathrm{F}}+C\bigl\|
\Lambda_{1}^{\ast}-\wt{\Lambda}_{1}\bigr\|
_{\mathrm{F}}.
\end{eqnarray*}
Here, the last inequality is due to
%
\begin{equation}
\label{eq:lambda-perturb} \| \wt{\Lambda}%
_{1}-\Lambda_{1}
\| _{\mathrm{F}}\leq\bigl\| \Sigma_{x}^{1/2}%
\bigl(
\wt{U}_{1}\wt{\Lambda}_{1}\wt{V}_{1}^{\prime}-U_{1}
\Lambda _{1}V_{1}^{\prime}\bigr)
\Sigma_{y}^{1/2}\bigr\| _{\mathrm{F}},
\end{equation}
a consequence of Lemma~\ref{lem:CCA2} and the  Hoffman--Wielandt inequality
\cite{tao12}, page 63.

We now control each of the three terms on the rightmost-hand side of
the second last display.
First, the bound we derived for (\ref{eq:zhaoren}),
up to a constant multiplier,
$\|\wt{U}_{1}\wt{\Lambda}_{1}\wt{V}_{1}^{\prime}-{U}_{1}\Lambda
_{1}{V}_{1}^{\prime}\| _{\mathrm{F}}$ is upper bounded by the
right-hand side of \eqref{eq:SparseApproClaim2}.
Next, the bound for $\|\Lambda_{1}^{\ast}-\wt{\Lambda}_{1}\|
_{\mathrm
{F}}$ has been shown in (\ref{eq:boundsLamdaTilda}).
Last but not least,
applying (\ref{eq:MPprincipal}) in Theorem~\ref{thm:sintheta}, we obtain
\begin{eqnarray*}
&&\bigl \|U_{1}^{\ast}\Lambda_{1}^{\ast}V_{1}^{\ast\prime}-
\wt {U}_{1}\wt {\Lambda}_{1}\wt{V}_{1}^{\prime}
\bigr\| _{\mathrm{F}}
\\
&&\qquad \leq\frac{C ( \wt{\lambda}_{1}+{\lambda}_{1}^*) }{\wt{\lambda
}_{r}-{\lambda}^*_{r+1}}\bigl\| (\Sigma_{xS_{u}S_{u}} )%
^{-1/2}
\Sigma_{xyS_{u}S_{v}} (\Sigma_{yS_{v}S_{v}} )^{-1/2}-\Xi \bigr\|
_{\mathrm{F}} 
\leq C \sqrt{\frac{q}{2-q} } \lambda
\eps_n %
,
\end{eqnarray*}
where the last inequality is due to (\ref{eq:usedlater4}), (\ref
{eq:boundsLamdaTilda}), (\ref{eq:lambda-lower}) and Lemma~\ref{lem:CCA2}.
The proof is completed by assembling the bounds for the three terms.
\end{pf*}

\subsection{Proof of Lemma \texorpdfstring{\protect\ref{lem:bias}}{4}}
\label{sec:bias-pf}

In this proof, we need the following technical result,
which is a direct consequence of Lemma~3 in \cite{supp2} by
applying union bound. Remember the notation $T_u$ and $T_v$ defined in
(\ref{eq:Tu-Tv}).

\begin{lemma} \label{lem:covdeviation45}
Assume
$\frac{1}{n}(k_{q}^{u}\log(ep/k_{q}^{u})+k_{q}^{v}\log
(em/k_{q}^{v}))<C_{1}$ for some constant $C_{1}>0$. For any constant
$C'>0$, there exists some constant $C>0$ only depending on $M,C_1$ and
$C'$, such that
\begin{eqnarray*}
\bigl\| \widehat{\Sigma}_{xT_{u}T_{u}}-\Sigma_{xT_{u}T_{u}}\bigr\|
_{\mathrm{op}}^{2} &\leq&\frac{C}{n}\bigl(k_{q}^{u}
\log\bigl(ep/k_{q}^{u}\bigr)\bigr),
\\
\bigl\| \widehat{\Sigma}_{yT_{v}T_{v}}-\Sigma_{yT_{v}T_{v}}\bigr\|
_{\mathrm{op}}^{2} &\leq&\frac{C}{n}\bigl(k_{q}^{v}
\log\bigl(em/k_{q}^{v}\bigr)\bigr),
\end{eqnarray*}
with probability at least $1-\exp(-C^{\prime}k_{q}^{u}\log
(ep/k_{q}^{u}))-\exp(-C^{\prime}k_{q}^{v}\log(em/k_{q}^{v}))$.
\end{lemma}

In addition, we need the following result.

\begin{lemma}[(Stewart and Sun \cite{Stewart90}, Theorem II.4.11)]
\label{lem:subspacedist} For any matrices $%
A,B $ with $A^{\prime}A=B^{\prime}B=I$, we have
\[
\inf_{W}\| A-BW\| _{\mathrm{F}}\leq\bigl\|
AA^{\prime
}-BB^{\prime}\bigr\| _{\mathrm{F}}.
\]
\end{lemma}

We first bound $ \langle\Sigma_xU_2\Lambda_2V_2'\Sigma_y, \wh
{U}_1\wh{V}_1'  \rangle$. By the definition of trace product,
we have
\begin{eqnarray*}
\bigl\langle\Sigma_xU_2\Lambda_2V_2'
\Sigma_y, \wh{U}_1\wh{V}_1'
\bigr\rangle&=& \bigl\langle\Lambda_2 V_2'
\Sigma_y\wh{V}_1', U_2'
\Sigma_x\wh{U}_1 \bigr\rangle
\\
&\leq& \bigl\|\Lambda_2 V_2'
\Sigma_y\wh{V}_1'\bigr\|_{\mathrm{F}}
\bigl\|U_2'\Sigma _x\wh{U}_1
\bigr\|_{\mathrm{F}}
\\
&\leq& \lambda_{r+1} \bigl\|V_2'
\Sigma_y\wh{V}_1'\bigr\|_{\mathrm{F}}
\bigl\|U_2'\Sigma _x\wh{U}_1
\bigr\|_{\mathrm{F}}.
\end{eqnarray*}
Define the SVD of matrices $U_1$ and $\wh{U}_1$ to be
\[
U_1=\Theta R H', \qquad\wh{U}_1=\wh{\Theta}
\wh{R}\wh{H}'.
\]
For any matrix $W$, we have
\begin{eqnarray*}
\bigl\|\wh{U}_1'\Sigma_xU_2
\bigr\|_{\mathrm{F}} &=& \bigl\|\bigl(\wh{U}_1-U_1HR^{-1}W
\wh {R}\wh{H}'\bigr)'\Sigma_xU_2
\bigr\|_{\mathrm{F}}
\\
&\leq& C\bigl\|\wh{U}_1-U_1HR^{-1}W\wh{R}
\wh{H}'\bigr\|_{\mathrm{F}}
\\
&\leq& C\|\wh{R}\|_{\mathrm {op}}\|\wh{\Theta}-\Theta W\|_{\mathrm{F}},
\end{eqnarray*}
where $\|\wh{R}\|_{\mathrm {op}}\leq\|\wh{U}_1\|_{\mathrm {op}}\leq\|(\wh
{\Sigma}_{xT_uT_u})^{-1/2}\|_{\mathrm{ op}}\|(\wh{\Sigma
}_{xT_uT_u})^{1/2}\wh{U}_{1T_u*}\|_{\mathrm{ op}}\leq C$ with probability
at least $1-\exp(-C^{\prime
}k_{q}^{u}\log(ep/k_{q}^{u}))-\exp(-C^{\prime}k_{q}^{v}\log
(em/k_{q}^{v}))$ by Lem\-ma~\ref{lem:covdeviation45}.
Hence, by Lem\-ma~\ref{lem:subspacedist}, we have
%
\begin{equation}
\bigl\|\wh{U}_1'\Sigma_xU_2
\bigr\|_{\mathrm{F}}\leq C\inf_W\|\wh{\Theta}-\Theta W
\|_{\mathrm{F}}\leq C\bigl\|\wh{\Theta}\wh{\Theta}'-\Theta
\Theta'\bigr\|_{\mathrm{F}}. \label{eq:fix1}
\end{equation}
We note that both $\widehat{\Theta}\widehat{\Theta}%
^{\prime}$ and $\Theta\Theta^{\prime}$ are the projection matrices of
the left singular spaces of $\widehat{U}_{1}\widehat{V}_{1}^{\prime}$
and $%
U_{1}V_{1}^{\prime}$, respectively, and the eigengap is at constant level
since the $r$th singular value of $U_{1}V_{1}^{\prime}$ is bounded below
by some constant and the $(r+1)$th singular value of $\widehat{U}_{1}%
\widehat{V}_{1}^{\prime}$ is zero.
Then a direct consequence of Wedin's sin-theta theorem \cite{Wedin72} gives
%
\begin{equation}
\bigl\|\wh{\Theta}\wh{\Theta}'-\Theta\Theta'
\bigr\|_{\mathrm{F}}\leq C\bigl\|\wh {U}_1\wh{V}_1'-U_1V_1'
\bigr\|_{\mathrm{F}}.\label{eq:fix2}
\end{equation}
Combining (\ref{eq:fix1}) and (\ref{eq:fix2}), we have $\|\wh
{U}_1'\Sigma_xU_2\|_{\mathrm{F}}\leq C_1\|\wh{U}_1\wh{V}_1'-U_1V_1'\|
_{\mathrm{F}}$. The same
argument also implies $\|V_2'\Sigma_y\wh{V}_1'\|_{\mathrm{F}}\leq C_1\|
\wh{U}_1\wh{V}_1'-U_1V_1'\|_{\mathrm{F}}$. Therefore,
\[
\bigl\llvert \bigl\langle\Sigma_xU_2
\Lambda_2V_2'\Sigma_y,
\wh{U}_1\wh {V}_1' \bigr\rangle\bigr
\rrvert \leq C_2\lambda_{r+1}\bigl\|\wh{U}_1
\wh{V}_1'-U_1V_1'
\bigr\|_{\mathrm{F}}^2.
\]
Using a similar argument, we also obtain
\[
\bigl\llvert \bigl\langle\Sigma_xU_2
\Lambda_2V_2'\Sigma_y, \wh
{U}_1^*\bigl(\wh{V}_1^*\bigr)' \bigr
\rangle\bigr\rrvert \leq C_2\lambda_{r+1}\bigl\|
\wh{U}^*_1\bigl(\wh{V}_1^*\bigr)'-U_1V_1'
\bigr\|_{\mathrm{F}}^2.
\]
By the  triangle inequality, we complete the proof.

\subsection{Proof of Lemma \texorpdfstring{\protect\ref{lem:excessloss1}}{5}}
\label{sec:excessloss1-pf}
Define
\[
W=%
\left[\matrix{ 0 & \widehat{U}_{1}^{\ast}
\widehat{V}_{1}^{\ast^{\prime}}-\widehat {U}_{1}%
\widehat{V}_{1}^{\prime}
\vspace*{2pt}\cr
\bigl(\widehat{U}_{1}^{\ast}
\widehat{V}_{1}^{\ast^{\prime}}-\widehat {U}_{1}%
\widehat{V}_{1}^{\prime}\bigr)^{\prime} & 0%
}
\right].
\]
Then simple algebra leads to 
%
\begin{equation}
\bigl\langle\Sigma_{xy}-\widehat{\Sigma}_{xy},\widehat
{U}_{1}^{\ast} 
\widehat{V}_{1}^{\ast^{\prime}}-
\widehat{U}_{1}\widehat {V}_{1}^{\prime
} \bigr
\rangle=\tfrac{1}{2} \langle\Sigma-\widehat{\Sigma}%
,W \rangle.
\label{eq:excessexpand}
\end{equation}
In the rest of the proof,
we bound $\langle\Sigma-\wh{\Sigma}, W \rangle$ by using Lemma~\ref{lem:EP}.

Notice that the matrix $\widehat{U}_{1}^{\ast}\widehat{V}_{1}^{\ast
^{\prime}}-\widehat{U}_{1}\widehat{V}_{1}^{\prime}$ has nonzero
rows with
indices in $T_{u}=S_{u}\cup\widehat{S}_{u}$ and nonzero columns with
indices in $T_{v}=S_{v}\cup\widehat{S}_{v}$. Hence, the enlarged
matrix $W$
has nonzero rows and columns with indices in $T\times T$, where
\[
T=T_{u} \cup(T_{v}+p)
\]
with $T_v +p = \{j+p: j\in T_v \}$.
The cardinality of $T$ is $|T| = |T_{u}|+|T_{v}| \leq
2(k_{q}^{u}+k_{q}^{v})$.
Thus, we can rewrite (\ref{eq:excessexpand}) as
\begin{eqnarray*}
&&\bigl\langle\Sigma_{xy}-\widehat{\Sigma}_{xy},\widehat
{U}_{1}^{\ast}%
\widehat{V}_{1}^{\ast^{\prime}}-
\widehat{U}_{1}\widehat {V}_{1}^{\prime
} \bigr
\rangle\\
&&\qquad= \tfrac{1}{2} \langle\Sigma-\wh{\Sigma}, W \rangle
\\
&&\qquad=\tfrac{1}{2} \langle\Sigma_{TT}-\widehat{\Sigma
}_{TT},W_{TT}\rangle
\\
&&\qquad=\tfrac{1}{2} \bigl\langle I_{|T|}-\Sigma_{TT}^{-1/2}
\widehat {\Sigma}%
_{TT}\Sigma_{TT}^{-1/2},
\Sigma_{TT}^{1/2}W_{TT}\Sigma _{TT}^{1/2}
\bigr\rangle
\\
&&\qquad=\tfrac{1}{2}\bigl\| \Sigma_{TT}^{1/2}W_{TT}
\Sigma _{TT}^{1/2}\bigr\| _{\mathrm{F}} \bigl\langle
I_{|T|}-\Sigma_{TT}^{-1/2}%
\widehat{
\Sigma}_{TT}\Sigma_{TT}^{-1/2}, K^T
\bigr\rangle,
\end{eqnarray*}
where $K^T = \frac{\Sigma
_{TT}^{1/2}W_{TT}\Sigma_{TT}^{1/2}}{\| \Sigma
_{TT}^{1/2}W_{TT}\Sigma_{TT}^{1/2}\| _{\mathrm{F}}}$.
Note that
%
\[
\tfrac{1}{2}\bigl\| \Sigma_{TT}^{1/2}W_{TT}
\Sigma_{TT}^{1/2}\bigr\| _{\mathrm{F}}\leq C\| W_{TT}
\| _{\mathrm{F}}=C\| W\| _{\mathrm{F}}=\sqrt{2}C\bigl\| \widehat{U}_{1}^{\ast}%
\widehat{V}_{1}^{\ast^{\prime}}-\widehat{U}_{1}\widehat
{V}_{1}^{\prime
}\bigr\| _{\mathrm{F}}.
\]
To obtain the desired bound,
it suffices to show that
%
\begin{equation}
\bigl\llvert \bigl\langle I_{|T|}-\Sigma_{TT}^{-1/2}
\widehat{\Sigma}%
_{TT}\Sigma_{TT}^{-1/2},
K^T \bigr\rangle\bigr\rrvert \label
{eq:useEPtobound}
\end{equation}
is upper bounded by $C\lambda\eps_n$ with high probability.

To this end, we note that
$T_{u} = S_u \cup\wh{S}_u$ has at most ${{p\choose k_{q}^{u}}}$
different possible configurations
since $S_u$ is deterministic and $\wh{S}_u$ is a random set of size $k_q^u$.
For the same reason, $T_{v}$
has at most ${{m\choose k_{q}^{v}}}$ different possible configurations.
Therefore, the subset $T$ has at
most $N={{p\choose k_{q}^{u}}}{{m\choose k_{q}^{v}}}$ different
possible configurations, which can be
listed as $T_{1},T_{2},\ldots,T_{N}$.
Let
\[
K^{T_j} = \frac{\Sigma_{T_{j}T_{j}}^{1/2}W_{T_{j}T_{j}}\Sigma
_{T_{j}T_{j}}^{1/2}}{\| \Sigma
_{T_{j}T_{j}}^{1/2}W_{T_{j}T_{j}}\Sigma_{T_{j}T_{j}}^{1/2}\| _{%
\mathrm{F}}}
\]
for all $j\in[N]$.
Since each $W_{T_{j}T_{j}}$ is of rank at most $2r$, so are the $K^{T_j}$'s.
Therefore,
\begin{eqnarray*}
\bigl|(\ref{eq:useEPtobound})\bigr| &\leq&\max_{1\leq j\leq N}\bigl\llvert \bigl\langle
I_{|T_{j}|}-\Sigma_{T_{j}T_{j}}^{-1/2}\widehat{\Sigma
}_{T_{j}T_{j}}\Sigma _{T_{j}T_{j}}^{-1/2}, K^{T_j} \bigr
\rangle\bigr\rrvert
\\
&\leq&\max_{1\leq j\leq N}\sup_{\| K\| _{\mathrm{F}%
}\leq1, \operatorname{rank}(K)\leq2r}\bigl\llvert
\bigl\langle I_{|T_{j}|}-\Sigma_{T_{j}T_{j}}^{-1/2}\widehat{
\Sigma }_{T_{j}T_{j}}\Sigma _{T_{j}T_{j}}^{-1/2},K \bigr\rangle\bigr
\rrvert .
\end{eqnarray*}
Then the union bound leads to
%
\begin{eqnarray}\label{eq:EPused1}
&&\mathbb{P}_\Sigma\bigl(\bigl|(\ref{eq:useEPtobound})\bigr|>t\bigr)
\nonumber
\\
&&\qquad \leq\sum_{j=1}^{N}\mathbb{P} \Bigl(\sup
_{\| K\| _{%
\mathrm{F}}\leq1,\operatorname{rank}(K)\leq2r}\bigl\llvert \bigl\langle I_{|T_{j}|}-
\Sigma_{T_{j}T_{j}}^{-1/2}\widehat{\Sigma }_{T_{j}T_{j}}\Sigma
_{T_{j}T_{j}}^{-1/2},K \bigr\rangle\bigr\rrvert >t \Bigr)
\\
&&\qquad\leq\sum_{j=1}^{N}\exp
\bigl(C^{\prime}r|T_{j}|-Cn\bigl(t\wedge t^{2}\bigr)
\bigr)
\nonumber\\
&&\qquad\leq{\pmatrix{p\cr k_{q}^{u}}} {\pmatrix{m\cr
k_{q}^{v}}}\exp\bigl(%
C_{1}r
\bigl(k_{q}^{u}+k_{q}^{v}\bigr)-Cn
\bigl(t\wedge t^{2}\bigr)\bigr)
\nonumber
\\
&&\qquad\leq\exp\biggl( C_{1}r\bigl(k_{q}^{u}+k_{q}^{v}
\bigr)+k_{q}^{u}\log \frac
{ep}{k_{q}^{u}}+k_{q}^{v}
\log\frac{em}{k_{q}^{v}}-Cn\bigl(t\wedge t^{2}\bigr) \biggr),
\nonumber
\end{eqnarray}
where  inequality (\ref{eq:EPused1}) is due to Lemma~\ref{lem:EP}.
We complete the proof by choosing
$t^{2}= C_2 \lambda^2 \eps_n^2
$ in the last display for some sufficiently large constant $C_2>0$,
which, by condition (\ref{eq:ass1}), is bounded.

\subsection{Proof of Lemma \texorpdfstring{\protect\ref{lem:excessloss2}}{6}}
\label{sec:excessloss2-pf}

First, we apply a telescoping expansion to the quantity of interest as
%
\begin{eqnarray}
&& \bigl\langle\widehat{\Sigma}_{x}\widehat{U}_{1}^{\ast}
\Lambda _{1}%
\widehat{V}_{1}^{\ast\prime}
\widehat{\Sigma}_{y}-{\Sigma}_{x}{U}%
_{1}
\Lambda_{1}{V}_{1}^{\prime}{\Sigma}_{y},
\widehat{U}_{1}^{\ast}%
\widehat{V}_{1}^{\ast^{\prime}}-
\widehat{U}_{1}\widehat {V}_{1}^{\prime
} \bigr
\rangle
\nonumber
\\
&&\qquad= \bigl\langle\Sigma_{x}\widehat{U}_{1}^{\ast}
\Lambda_{1}\widehat {V}%
_{1}^{\ast\prime}
\Sigma_{y}-\Sigma_{x}U_{1}^{\ast}
\Lambda _{1}V_{1}^{\ast\prime}\Sigma_{y},
\widehat{U}_{1}^{\ast}\widehat{V} 
_{1}^{\ast^{\prime}}-
\widehat{U}_{1}\widehat{V}_{1}^{\prime} \bigr
\rangle%
\label{eq:excessloss2.3}
\\
&&\qquad\quad{}+ \bigl\langle\Sigma_{x}U_{1}^{\ast}
\Lambda_{1}V_{1}^{\ast\prime
}\Sigma_{y}-{
\Sigma}_{x}{U}_{1}\Lambda_{1}{V}_{1}^{\prime}{
\Sigma }_{y},%
\widehat{U}_{1}^{\ast}
\widehat{V}_{1}^{\ast^{\prime}}-\widehat{U}_{1} 
\widehat{V}_{1}^{\prime} \bigr\rangle\label{eq:excessloss2.2}
\\
&&\qquad\quad{}+ \bigl\langle\widehat{\Sigma}_{x}\widehat{U}_{1}^{\ast}
\Lambda_{1} 
\widehat{V}_{1}^{\ast\prime}
\widehat{\Sigma}_{y}-\Sigma _{x}\widehat
{U}%
_{1}^{\ast}\Lambda_{1}
\widehat{V}_{1}^{\ast\prime}\Sigma _{y},
\widehat{U}%
_{1}^{\ast}\widehat{V}_{1}^{\ast^{\prime}}-
\widehat{U}_{1}\widehat{V} 
_{1}^{\prime}
\bigr\rangle. \label{eq:excessloss2.1}
\end{eqnarray}
In what follows, we bound each of the terms in \eqref
{eq:excessloss2.3}--\eqref{eq:excessloss2.1} in order.


\textit{$1^\circ$ Bound for} \eqref{eq:excessloss2.3}.
Applying (\ref{eq:OracleLossClaim2}) in Lemma~\ref{lem:oracleloss},
we obtain
that with probability at least $1-\exp(-C^{\prime}(k_{q}^{u}+\log
(ep/k_{q}^{u})))-\exp(-C^{\prime}(k_{q}^{v}+\log
(em/k_{q}^{v})))$,
\begin{eqnarray*}
\bigl\llvert (\ref{eq:excessloss2.3})\bigr\rrvert &\leq&C\bigl\| \widehat{U
}_{1}^{\ast}\Lambda_{1}
\widehat{V}_{1}^{\ast\prime}-U_{1}^{\ast
}\Lambda
_{1}V_{1}^{\ast\prime}\bigr\| _{\mathrm{F}}\bigl\|
\widehat{U}%
_{1}^{\ast}\widehat{V}_{1}^{\ast^{\prime}}-
\widehat{U}_{1}\widehat{V} 
_{1}^{\prime}
\bigr\| _{\mathrm{F}}
\\
&\leq& 
C\sqrt{
\frac{q}{2-q}}\lambda\eps_n 
\bigl\|
\widehat{U}_{1}^{\ast}\widehat{V}_{1}^{\ast^{\prime}}-
\widehat{ 
U}_{1}\widehat{V}_{1}^{\prime}
\bigr\| _{\mathrm{F}}.
\end{eqnarray*}

\textit{$2^\circ$ Bound for} \eqref{eq:excessloss2.2}.
Applying (\ref{eq:SparseApproClaim2}) in Lemma~\ref
{lem:sparseapproxerror}, we obtain
\begin{eqnarray*}
\bigl\llvert (\ref{eq:excessloss2.2})\bigr\rrvert &\leq&C\bigl\| U_{1}^{\ast}
\Lambda_{1}V_{1}^{\ast\prime}-{U}_{1}
\Lambda_{1}{V}%
_{1}^{\prime}\bigr\|
_{\mathrm{F}}\bigl\| \widehat{U}_{1}^{\ast}%
\widehat{V}_{1}^{\ast^{\prime}}-\widehat{U}_{1}\widehat
{V}_{1}^{\prime
}\bigr\| _{\mathrm{F}}
\\
&\leq& C\sqrt{\frac{q}{2-q}}\lambda\eps_n \bigl\|
\widehat{U}_{1}^{\ast}\widehat{V}_{1}^{\ast^{\prime}}-
\widehat{ 
U}_{1}\widehat{V}_{1}^{\prime}
\bigr\| _{\mathrm{F}}. 
\end{eqnarray*}
%

\textit{$3^\circ$ Bound for} \eqref{eq:excessloss2.1}.
We turn to bound (\ref{eq:excessloss2.1}) based on a strategy similar
to that used in proving Lemma~\ref{lem:excessloss1}.
First, we write
it in a form for which we could apply Lemma~\ref{lem:EP}.
Recall the random sets $T_u$ and $T_v$ defined in \eqref{eq:Tu-Tv}.
Then for
\begin{eqnarray*}
H_{x}^{T_u} & =&(\Sigma_{xT_{u}T_{u}})^{1/2}
\bigl(\widehat{U}_{1T_{u}\ast
}^{\ast}\bigl(\widehat{V}_{1T_{v}\ast}^{\ast}
\bigr)^{\prime}-\widehat{U}%
_{1T_{u}\ast}(\widehat{V}_{1T_{v}\ast})^{\prime}
\bigr)
\\
&&{} \times\widehat{\Sigma}%
_{yT_{v}T_{v}}\widehat{V}_{1T_{v}\ast}^{\ast}
\Lambda_{1}\bigl(\widehat {U}%
_{1T_{u}\ast}^{\ast}
\bigr)^{\prime}(\Sigma_{xT_{u}T_{u}})^{1/2},
\\
H_{y}^{T_v} & = &(\Sigma_{yT_{v}T_{v}})^{1/2}
\widehat{V}_{1T_{v}\ast
}^{\ast
}\Lambda_{1}\bigl(
\widehat{U}_{1T_{u}\ast}^{\ast}\bigr)^{\prime}
\\
&&{} \times\Sigma _{xT_{u}T_{u}}\bigl(\widehat{U}_{1T_{u}\ast}^{\ast}
\bigl(\widehat {V}_{1T_{v}\ast
}^{\ast}\bigr)^{\prime}-
\widehat{U}_{1T_{u}\ast}(\widehat {V}_{1T_{v}\ast
})^{\prime}\bigr) (
\Sigma_{yT_{v}T_{v}})^{1/2},
\end{eqnarray*}
and $\overline{H}_x^{T_u} = {H_{x}^{T_u}}/{\|H_{x}^{T_u}\| _{\mathrm
{F}}}$, $\overline{H}_y^{T_v} = {H_{y}^{T_v}}/{\|H_{y}^{T_v}\|
_{\mathrm
{F}}}$, we have
\begin{eqnarray*}
\label{eq:Excess2FirstDec}
&&\bigl\llvert (\ref{eq:excessloss2.1})\bigr\rrvert
\nonumber
\\
&&\qquad= \bigl\llvert \bigl\langle\widehat{\Sigma}_{x}-
\Sigma_{x},\bigl(\widehat {U}_{1}^{\ast
}
\widehat{V}_{1}^{\ast\prime}-\widehat{U}_{1}\widehat
{V}_{1}{^{\prime
}}%
\bigr)\widehat{
\Sigma}_{y}\widehat{V}_{1}^{\ast}
\Lambda_{1}\widehat {U}_{1}^{\ast}{%
^{\prime}}
\bigr\rangle
\nonumber
\\
&&\qquad\quad{}  + \bigl\langle\widehat{\Sigma}_{y}-\Sigma
_{y},\widehat{V}_{1}^{\ast}\Lambda_{1}
\widehat{U}_{1}^{\ast
}{^{\prime
}}\Sigma_{x}%
\bigl(\widehat{U}_{1}^{\ast}\widehat{V}_{1}^{\ast\prime}-
\widehat {U}_{1}%
\widehat{V}_{1}{^{\prime}}
\bigr) \bigr\rangle\bigr\rrvert
\nonumber
\\
&&\qquad\leq \bigl\llvert \bigl\langle\widehat{\Sigma}_{xT_{u}T_{u}}-\Sigma
_{xT_{u}T_{u}},\bigl(%
\widehat{U}_{1T_{u}\ast}^{\ast}
\bigl(\widehat{V}_{1T_{v}\ast}^{\ast
}\bigr)^{\prime
}\\
&&\qquad\quad{}-
\widehat{U}_{1T_{u}\ast}(\widehat{V}_{1T_{v}\ast})^{\prime}
\bigr)%
\widehat{\Sigma}_{yT_{v}T_{v}}\widehat{V}_{1T_{v}\ast}^{\ast
}
\Lambda _{1}\bigl(%
\widehat{U}_{1T_{u}\ast}^{\ast}
\bigr)^{\prime} \bigr\rangle\bigr\rrvert
\nonumber
\\
&&\qquad\quad{}+\bigl\llvert \bigl\langle\widehat{\Sigma}_{yT_{v}T_{v}}\\
&&\qquad\quad{}-\Sigma
_{yT_{v}T_{v}},%
\widehat{V}_{1T_{v}\ast}^{\ast}
\Lambda_{1}\bigl(\widehat{U}_{1T_{u}\ast
}^{\ast}
\bigr)^{\prime}\Sigma_{xT_{u}T_{u}}\bigl(\widehat{U}_{1T_{u}\ast
}^{\ast}
\bigl(\widehat{V}_{1T_{v}\ast}^{\ast}\bigr)^{\prime}-
\widehat{U}%
_{1T_{u}\ast}(\widehat{V}_{1T_{v}\ast})^{\prime}
\bigr) \bigr\rangle \bigr\rrvert
\nonumber
\\
&&\qquad=\bigl\| H_{x}^{T_u}\bigr\| _{\mathrm{F}}\bigl\llvert \bigl
\langle(\Sigma _{xT_{u}T_{u}})^{-1/2}\widehat{\Sigma}_{xT_{u}T_{u}}(
\Sigma _{xT_{u}T_{u}})^{-1/2}-I_{|T_{u}|}, \overline{H}_x^{T_u}
\bigr\rangle \bigr\rrvert
\nonumber
\\
&&\qquad\quad{}+\bigl\| H_{y}^{T_v}\bigr\| _{\mathrm{F}} \bigl\llvert \bigl
\langle(\Sigma _{yT_{v}T_{v}})^{-1/2}\widehat{\Sigma}_{yT_{v}T_{v}}(
\Sigma _{yT_{v}T_{v}})^{-1/2}-I_{|T_{v}|}, \overline{H}_y^{T_v}
\bigr\rangle \bigr\rrvert .
\nonumber
\end{eqnarray*}

We now bound each term on the rightmost side.
Applying Lemma~\ref{lem:EP} with union bound and then following a
similar analysis
to that leading to (\ref{eq:useEPtobound}) but with $T$ replaced by
$T_{u}$ and $T_{v}$, we
obtain that
%
\begin{eqnarray}\label{eq:excess2.1.1}
%
 &&\bigl\llvert \bigl\langle(\Sigma
_{xT_{u}T_{u}})^{-1/2}\widehat{\Sigma}_{xT_{u}T_{u}}(\Sigma
_{xT_{u}T_{u}})^{-1/2}-I_{|T_{u}|}, \overline{H}_x^{T_u}
\bigr\rangle \bigr\rrvert \nonumber\\
&&\qquad\leq C\sqrt{\frac{k_{q}^{u}%
}{n}\biggl( r+\log
\frac{ep}{k_{q}^{u}} \biggr)},
\nonumber
\\[-8pt]
\\[-8pt]
\nonumber
&&\bigl\llvert \bigl\langle(\Sigma _{yT_{v}T_{v}})^{-1/2}\widehat{
\Sigma}_{yT_{v}T_{v}}(\Sigma _{yT_{v}T_{v}})^{-1/2}-I_{|T_{v}|},
\overline{H}_y^{T_v} \bigr\rangle \bigr\rrvert\\
&&\qquad \leq C\sqrt
{\frac{k_{q}^{v}}{n}\biggl( r+\log\frac{em}{k_{q}^{v}} \biggr)}\nonumber
\end{eqnarray}
with probability at least $1-\exp(-C'k_{q}^{u}(r+\log(ep/k_{q}^{u})%
))$ and $1-\exp(-C'k_{q}^{v}(r+\log(em/k_{q}^{v})))$,
respectively.


To bound $\| H_{x}^{T_u}\| _{\mathrm{F}}$ and $\|
H_{y}^{T_v}\| _{\mathrm{F}}$,
we note that it follows from Lemma~\ref{lem:covdeviation45} that all eigenvalues of $\widehat{\Sigma}%
_{xT_{u}T_{u}}$ and $\widehat{\Sigma}_{yT_{v}T_{v}}$ are bounded from
below and
above by some universal positive constants with probability at least
$1-\exp(%
-C^{\prime}k_{q}^{u}\log(ep/k_{q}^{u}))-\exp(-C^{\prime
}k_{q}^{v}\log(em/k_{q}^{v}))$ under assumption (\ref{eq:ass1}).
Thus, with the same probability we
have
%
\begin{eqnarray}
\label{eq:excess2.1.2} \bigl\| H_{x}^{T_u}\bigr\| _{\mathrm{F}} &\leq& C
\lambda\bigl\| \widehat{U}_{1}^{\ast}\widehat{V}_{1}^{\ast\prime}-
\widehat{U}_{1} 
\widehat{V}_{1}{^{\prime}}
\bigr\| _{\mathrm{F}}\bigl\| \widehat{\Sigma }_{yT_{v}T_{v}}^{1/2}
\widehat{V}_{1T_{v}\ast}^{\ast}\bigr\| _{\mathrm
{op}}
\nonumber
\\
&&{}\times \bigl\| \widehat{\Sigma}_{yT_{v}T_{v}}^{1/2}\bigr\| _{%
\mathrm{op}}\bigl\|
\widehat{\Sigma}_{xT_{u}T_{u}}^{1/2}\widehat{U}%
_{1T_{u}\ast}^{\ast}
\bigr\| _{\mathrm{op}}\bigl\| \widehat{\Sigma}%
_{xT_{u}T_{u}}^{-1/2}\bigr\|
_{\mathrm{op}}
\\
&\leq&C_1\lambda\bigl\| \widehat{U}_{1}^{\ast}
\widehat{V}%
_{1}^{\ast\prime}-\widehat{U}_{1}
\widehat{V}_{1}{^{\prime}}\bigr\| _{%
\mathrm{F}}\nonumber
\end{eqnarray}
and
%
\begin{eqnarray}
\label{eq:excess2.1.3}\bigl \| H_{y}^{T_v}\bigr\| _{\mathrm{F}} &\leq& C
\lambda\bigl\| \widehat{U}_{1}^{\ast}\widehat{V}_{1}^{\ast\prime}-
\widehat {U}_{1}\widehat{V}_{1}{^{\prime}}\bigr\|
_{\mathrm{F}}\bigl\| \widehat{\Sigma }_{yT_{v}T_{v}}^{1/2}
\widehat{V}_{1T_{v}\ast}^{\ast}\bigr\| _{%
\mathrm{op}}
\nonumber
\\
&&{}\times \bigl\| \widehat{\Sigma}_{yT_{v}T_{v}}^{-1/2}\bigr\| _{%
\mathrm{op}}\bigl\|
\widehat{\Sigma}_{xT_{u}T_{u}}^{1/2}\widehat{U}%
_{1T_{u}\ast}^{\ast}
\bigr\| _{\mathrm{op}}\bigl\| \widehat{\Sigma}%
_{xT_{u}T_{u}}^{-1/2}\bigr\|
_{\mathrm{op}}
\\
&\leq&C_1\lambda\bigl\| \widehat{U}_{1}^{\ast}
\widehat{V}%
_{1}^{\ast\prime}-\widehat{U}_{1}
\widehat{V}_{1}{^{\prime}}\bigr\| _{%
\mathrm{F}}.\nonumber
\end{eqnarray}
Combining (\ref{eq:excess2.1.1}), (\ref{eq:excess2.1.2}) and (\ref%
{eq:excess2.1.3}), we obtain
%
\[
\bigl\llvert (\ref{eq:excessloss2.1})\bigr\rrvert \leq C\lambda^2 \eps
_n \bigl\| \widehat{U}_{1}^{\ast}
\widehat{V}_{1}^{\ast\prime}-\widehat{U}_{1}
\widehat{V}_{1}{^{\prime}}\bigr\| _{\mathrm{F}},
\]
with probability at
least $1-\exp(-C^{\prime}k_{q}^{u}\log(ep/k_{q}^{u}))-\exp(%
-C^{\prime}k_{q}^{v}\log(em/k_{q}^{v}))$.
Noting that $\lambda< 1$, this completes the proof.


\begin{supplement}[id=suppA]
\stitle{Supplement to ``Minimax estimation in sparse canonical
correlation analysis''}
\slink[doi]{10.1214/15-AOS1332SUPP} 
\sdatatype{.pdf}
\sfilename{aos1332\_supp.pdf}
\sdescription{The supplement \cite{supp2} contains an Appendix to the current paper
in which we prove
Theorems \ref{thm:lower-bd-q}--\ref{thm:sintheta} and Lemmas \ref
{lem:oracleloss} and \ref{lem:EP}--\ref{lem:lq3}.}
\end{supplement}

%



\printaddresses
\end{document}